\documentclass[10pt,reqno,twoside]{amsart}

\usepackage{amsfonts}
\usepackage{paralist}
  \usepackage{graphics}
  \usepackage{epsfig}
 \usepackage[colorlinks=true]{hyperref}
\usepackage{hyperref}

\numberwithin{equation}{section}
\usepackage{latexsym}
\usepackage{amsmath}
\usepackage{amssymb}
\usepackage{mathrsfs}
\usepackage{graphicx}
\usepackage{ifthen}
\usepackage[T1]{fontenc}
\usepackage[latin1]{inputenc}
\usepackage{color}
\newtheorem{Satz}{Theorem}[section]

\newtheorem{Def}[Satz]{Definition}

\newcommand{\oo}{\overline \Omega}

\newcommand{\po}{\partial\Omega}

\newcommand{\diag}{\text {diag}}

\newcommand{\diam}{\text {diam}}

\newcommand{\supp}{\text {supp}}

\begin{document}

\title[The Lazer-McKenna conjecture]
{The Lazer-McKenna conjecture for an anisotropic planar exponential nonlinearity with a singular source}


\maketitle

%

%
%
%
%


\centerline{\scshape Yibin Zhang
}
\medskip
{\footnotesize
 \centerline{College of Sciences, Nanjing Agricultural University, Nanjing 210095, China}
  \smallskip
 \centerline{Email:\,\,\,\,yibin10201029@njau.edu.cn}
}

\begin{abstract}
Given a bounded smooth domain  $\Omega$ in $\mathbb{R}^2$,
we study the following  anisotropic elliptic  problem
$$
\begin{cases}
-\nabla\big(a(x)\nabla \upsilon\big)=
a(x)\big[e^{\upsilon}-s\phi_1-4\pi\alpha\delta_q-h(x)\big]\,\,\,\,
\,\textrm{in}\,\,\,\,\,\Omega,\\[2mm]
\upsilon=0
\qquad\qquad\qquad\qquad\qquad
\qquad\qquad\qquad\qquad\quad
\textrm{on}\,\ \,\partial\Omega,
\end{cases}
$$
where $a(x)$ is a positive smooth function,
$s>0$ is a large  parameter,
$h\in C^{0,\gamma}(\overline{\Omega})$,
$q\in\Omega$, $\alpha\in(-1,+\infty)\setminus\mathbb{N}$,
$\delta_q$ denotes the Dirac measure with pole at point $q$
and $\phi_1$ is a positive first
eigenfunction of
the problem
$-\nabla\big(a(x)\nabla \phi\big)=\lambda a(x)\phi$
under Dirichlet boundary condition in $\Omega$.
We show that  if $q$ is  both
a local maximum point of $\phi_1$
and  an isolated local maximum
point of $a(x)\phi_1$,
this problem  has a family of solutions $\upsilon_s$ with  arbitrary $m$ bubbles
accumulating to  $q$
and the quantity $\int_{\Omega}a(x)e^{\upsilon_s}\rightarrow8\pi(m+1+\alpha)a(q)\phi_1(q)$
as $s\rightarrow+\infty$,  which give a positive answer
to the Lazer-McKenna conjecture for this case.
\\
\noindent{
2020 \it Mathematics Subject Classification.}\,\,
Primary 35B25, 35J25; Secondary 35B40.\\
\noindent{\it Keywords:}\,\,
Lazer-McKenna conjecture;
Anisotropic coefficient;
Exponential nonlinearity;
Singular source;
Lyapunov-Schmidt procedure.
\end{abstract}

\maketitle


\section{Introduction}
We consider the  anisotropic elliptic   problem
\begin{equation}\label{1.1}
\aligned
\left\{\aligned
&-\nabla\big(a(x)\nabla \upsilon\big)=
a(x)\big[e^{\upsilon}-s\phi_1-4\pi\alpha\delta_q-h(x)\big]\,\,\,\,
\,\textrm{in}\,\,\,\,\,\Omega,\\[2mm]
&\upsilon=0
\qquad\qquad\qquad\qquad\qquad
\qquad\qquad\qquad\qquad\quad\
\textrm{on}\,\ \,\partial\Omega,
\endaligned\right.
\endaligned
\end{equation}
where
 $\Omega$ is a bounded smooth domain in $\mathbb{R}^2$,
$s>0$ is a large  parameter, $q\in\Omega$, $\alpha\in(-1,+\infty)\setminus\mathbb{N}$,
$\delta_q$ denotes the Dirac measure with pole at point $q$,
$h\in C^{0,\gamma}(\overline{\Omega})$ is given,
$a(x)$ is a smooth function over $\overline{\Omega}$ satisfying
\begin{equation}\label{1.5}
\aligned
a_1\leq a(x)\leq a_2
\endaligned
\end{equation}
for some constants $0<a_1<a_2<+\infty$, $\phi_1$ is a positive eigenfunction of
the anisotropic Laplacian operator
\begin{equation}\label{1.6}
\aligned
-\Delta_a:=-\frac{1}{a(x)}\nabla\big(a(x)\nabla\cdot\big)=-\Delta-\nabla\log a(x)\nabla
\endaligned
\end{equation}
with Dirichlet boundary condition corresponding to the
first eigenvalue $\lambda_1$.
Clearly, if we set $\rho(x)=(-\Delta_a)^{-1}h$ in $H_0^1(\Omega)$ and
let $G(x,y)$ be the Green's function satisfying
\begin{equation}\label{1.2}
\left\{\aligned
&-\Delta_aG(x,y)=8\pi\delta_{y}(x),
\,\,\,\,\,\,\,\,
x\in\Omega,\\[2mm]
&G(x,y)=0,
\,\qquad\qquad\quad\,\,\quad\,
x\in\partial\Omega,
\endaligned\right.
\end{equation}
and $H(x,y)$ be its regular part given by
\begin{equation}\label{1.3}
\aligned
H(x,y)=G(x,y)-4\log\frac1{|x-y|},
\endaligned
\end{equation}
then equation (\ref{1.1}) is equivalent to solving for
$u=\upsilon+\frac{s}{\lambda_1}\phi_1+\frac{\alpha}{2}G(\,\cdot,\,q)+\rho$, the problem
\begin{equation}\label{1.4}
\aligned
\left\{\aligned
&-\Delta_a u=|x-q|^{2\alpha}k(x)
e^{-t\phi_1}e^u\,\,\,\,\,
\textrm{in}\,\,\,\,\,\Omega,\\[2mm]
&u=0\,\,\qquad\,\,\qquad\,\,
\quad\quad\ \,\qquad\qquad\,\textrm{on}
\,\,\,\partial\Omega,
\endaligned\right.
\endaligned
\end{equation}
where $k(x)=e^{-\rho(x)-\frac{\alpha}{2}H(x,q)}$ and
$t=s/\lambda_1$.
The point $q$ with Dirac measure in equation (\ref{1.1})
is called    vortex point or singular source.
The term $|\cdot|^{2\alpha}$ in problem (\ref{1.4}) is called
the Hardy weight if $-1<\alpha<0$, whereas the H\'{e}non weight  if $\alpha>0$.
In this paper, we are interested in the existence of solutions of problem (\ref{1.4})
(or (\ref{1.1}))
which exhibit the {\it multiple concentration behavior} around singular source $q$ when
the  parameter $t$ tends to infinity.


If $a(x)=1$ and $\alpha=0$, equation (\ref{1.1}) is  the classical
 Ambrosetti-Prodi type problem \cite{AP}  with  exponential nonlinearity
\begin{equation}\label{1.7}
\left\{\aligned
&-\Delta \upsilon=g(\upsilon)-s\phi_1-h(x)\,\,\,\,
\,\,\textrm{in}\,\,\,\,\,\Omega,\\[1mm]
&\upsilon=0\,\,\qquad\qquad\quad\qquad\,
\qquad\,\,\,\,\textrm{on}\,\,\,\,\po,
\endaligned\right.
\end{equation}
as $s\rightarrow+\infty$, where $\Omega$ is a bounded smooth domain in $\mathbb{R}^N$
($N\geq2$),
$h\in C^{0,\gamma}(\overline{\Omega})$ is given,
$\phi_1$ is a positive eigenfunction of $-\Delta$ with Dirichlet
boundary condition corresponding to the first eigenvalue $\lambda_1$,
and
$g:\mathbb{R}\rightarrow\mathbb{R}$ is a continuous function such that
$-\infty\leq\nu=\lim_{t\rightarrow-\infty}\frac{g(t)}{t}<
\lim_{t\rightarrow+\infty}\frac{g(t)}{t}=\overline{\nu}\leq+\infty$.
Here $\nu=-\infty$
and $\overline{\nu}=+\infty$ are allowed.
The condition that $(\nu, \overline{\nu})$ contains some eigenvalues of $-\Delta$
subject to Dirichlet boundary condition
has great influence on   the existence and multiplicity
of solutions for problem (\ref{1.7}),
which has been   widely considered  since the early 1970s
so that  many interesting results have been obtained
(see \cite{D1,DS1,LM,RS3,S} and references therein).
Moreover, in the early 1980s Lazer and McKenna \cite{LM} conjectured that the
number of solutions to problem (\ref{1.7}) is
 unbounded as $s\rightarrow+\infty$ when
$\nu<\lambda_1<\overline{\nu}=+\infty$ and $g(t)$ has
an appropriate growth at infinity.

In fact,
the Lazer-McKenna conjecture
holds true for problem (\ref{1.7}) with many different types of nonlinearities,
which has been founded
in \cite{DM,DHZ,YZ} for the exponential nonlinear case
$g(t)=e^t-4\pi\alpha\delta_q$, $-1<\alpha\not\in\mathbb{N}^*$ and $N\geq2$,
in \cite{MP1,MP2,MP3} for the asymptotically linear case
$g(t)=\overline{\nu} t_{+}-\nu t_{-}$ with $\overline{\nu}$ large enough,
in \cite{DY2,DY1} for the subcritical  case
$g(t)=t_{+}^p+\lambda t$, $\lambda<\lambda_1$, $1<p<\frac{N+2}{N-2}$ if $N\geq3$,
$1<p<+\infty$ if $N=2$,
in \cite{D,LYY1,LYY2,WY2} for the critical case $g(t)=t_{+}^{\frac{N+2}{N-2}}+\lambda t$,
$0<\lambda<\lambda_1$ and $N\geq6$, in \cite{DS} for the superlinear
nonhomogeneous case $g(t)=t_{+}^p+t_{-}^q$,
$1<q<p<\frac{N+2}{N-2}$ and $N\geq4$,
in \cite{BMP,DY0,DSS} for the
superlinear
homogeneous case $g(t)=|t|^p$, $1<p<\frac{N+2}{N-2}$ if $N\geq3$,
$1<p<+\infty$ if $N=2$,
 and in \cite{DY3,WY1,ADM,GL,LP} for
some other  cases and even some generalized versions,
where $t_{+}=\max\{t,0\}$ and $t_{-}=\max\{-t,0\}$.
 In particular,
 when  $g(t)=e^t$ and $N=2$,
del Pino and Mu\~{n}oz in \cite{DM} proved the Lazer-McKenna conjecture
by constructing  solutions of problem (\ref{1.7})
with the accumulation of
arbitrarily many
bubbles around any isolated local maximum point of $\phi_1$.
Recently,  in \cite{DHZ}, Dong, Hu and Zhang extended this result in \cite{DM}
to the case $g(t)=e^t-4\pi\alpha\delta_q$, $-1<\alpha\not\in\mathbb{N}^*$ and $N=2$,
and showed that if singular source  $q$ is an isolated local maximum point of $\phi_1$,
problem (\ref{1.7}) always admits a family of solutions with  arbitrarily  many bubbles accumulating
to $q$.

It is necessary to point out that the anisotropic
equation (\ref{1.1})$|_{\alpha=0}$ is a special case
of problem (\ref{1.7}) with $g(t)=e^t$ in higher dimension $N\geq3$.
Indeed, let   a standard $N$-dimensional torus be
$\mathbb{T}=\{(x',x_N)\in\mathbb{R}^N:\,(\|x'\|-1)^2+x_N^2\leq r_0^2\}$
with $x'=(x_1,\ldots,x_{N-1})$,
$\|x'\|=\sqrt{x_1^2+\ldots+x_{N-1}^2}$ and $0<r_0<1$.
If we look for some special solutions of problem (\ref{1.7}) with $g(t)=e^t$ in the
 axially symmetric
torus $\mathbb{T}$,
i.e. solutions $\upsilon$ of the form $\upsilon(x',x_N)=\upsilon(r,x_N)$
with $r=\|x'\|$, a simple calculation shows that problem (\ref{1.7}) with $g(t)=e^t$
 in higher dimension  $N\geq3$ is reduced to
$$
\aligned
\left\{\aligned
&-\nabla\big(r^{N-2}\nabla \upsilon\big)=
r^{N-2}\big[e^{\upsilon}-s\phi_1-h(x)\big]\ \,\,\,
\,\textrm{in}\,\,\,\,\,\Omega_{\mathbb{T}},\\[2mm]
&\upsilon=0
\,\,\ \quad\quad\quad\quad\quad\quad\quad\,\,
\quad\quad\quad\quad\quad\,\,\,\,\,\,
\quad\quad\,\,\textrm{on}\,\ \,\partial\Omega_{\mathbb{T}},
\endaligned\right.
\endaligned
$$
where $\Omega_{\mathbb{T}}=\{(r,x_N)\in\mathbb{R}^2:\,(r-1)^2+x_N^2<r_0^2\}$.
This is just the equation (\ref{1.1})$|_{\alpha=0}$  with
anisotropic coefficient $a(r,x_N)=r^{N-2}$.
In this direction,
 Yang and Zhang in \cite{YZ}
proved that the Lazer-McKenna conjecture is also true
for  problem (\ref{1.7}) with $g(t)=e^t$ in the higher-dimensional
domain with some rotational symmetries, by constructing solutions of
 the anisotropic
equation (\ref{1.1})$|_{\alpha=0}$ which exhibit multiple
concentration behavior around local maximum points of $a(x)\phi_1$
in the domain as $s\rightarrow+\infty$.

In the present paper, our goal  is
to give a positive answer to the Lazer-McKenna conjecture
for the singular case of the anisotropic equation
(\ref{1.1})
involving $\alpha\in(-1,+\infty)\setminus\mathbb{N}$,
 by trying to prove the existence of multiple clustered  blowup
  solutions of problem (\ref{1.4})
with Hardy-H\'{e}non weight $|\cdot-q|^{2\alpha}$
 in a constructive way.
As a result, we  find that
if
 singular source $q$
 is  both
a local maximum point of $\phi_1$
and  an isolated local maximum
point of $a(x)\phi_1$,
 problem (\ref{1.4}) (or (\ref{1.1})) always admits a
family of solutions with  arbitrarily many bubbles
accumulating to  $q$.
In particular, we recover and extend the results in \cite{DM,DHZ,YZ}.
This  can be stated as follows.

\vspace{1mm}
\vspace{1mm}
\vspace{1mm}
\vspace{1mm}

\noindent{\bf Theorem 1.1.} {\it
Let $\alpha\in(-1,+\infty)\setminus\mathbb{N}$ and assume\\
\indent {\upshape(A)}\,\, Singular source
 $q$ is  both
a local maximum point of $\phi_1$
and  an isolated local maximum
point of $a(x)\phi_1$.
\\
Then for any integer $m\geq1$, there exists $t_m>0$
such that for any $t>t_m$,  problem {\upshape(\ref{1.4})} has
a family of solutions $u_t$
satisfying
$$
\aligned
u_t(x)=\left[\log
\frac{1}{(\varepsilon_{0,t}^2\mu_{0,t}^2+|x-q|^{2(1+\alpha)})^2}
+(1+\alpha)H(x,q)\right]
+
\sum\limits_{i=1}^{m}\left[\,\log
\frac1{(\varepsilon_{i,t}^2\mu_{i,t}^2+|x-\xi_{i,t}|^2)^2}
+H(x,\xi_{i,t})
\,\right]+o(1),
\endaligned
$$
where $o(1)\rightarrow0$, as $t\rightarrow+\infty$, uniformly on
each compact subset of $\overline{\Omega}\setminus\{q,\,\xi_{1,t},\ldots,\xi_{m,t}\}$,
the parameters $\varepsilon_{0,t}$, $\varepsilon_{i,t}$,  $\mu_{0,t}$  and
$\mu_{i,t}$ satisfy
$$
\aligned
\varepsilon_{0,t}=e^{-\frac12t\phi_1(q)},
\,\,\ \quad\ \,\,
\varepsilon_{i,t}=e^{-\frac12t\phi_1(\xi_{i,t})},
\,\,\ \quad\ \,\,
\frac1C\leq\mu_{0,t}\leq Ct^{\frac{m(m+1+\alpha)^2a_2}{a_1}},
\,\,\ \quad\ \,\,
\frac1C\leq\mu_{i,t}\leq Ct^{\frac{(2m+\alpha)(m+1+\alpha)^2a_2}{2a_1}},
\endaligned
$$
for  some   $C>0$,
and $(\xi_{1,t},\ldots,\xi_{m,t})\in\Omega^m$ satisfies
$$
\aligned
\xi_{i,t}\rightarrow q\,\ \ \,\,\textrm{for all}\,\,\,i,\,\,
\ \ \ \,
\textrm{and}\,\,\ \ \,\,\,|\xi_{i,t}-\xi_{j,t}|>t^{-\frac{(m+1+\alpha)^2a_2}{2a_1}}\,\,\ \,\,\forall
\,\,\,
i\neq j.
\endaligned
$$
}

The corresponding result for problem (\ref{1.1}) can be stated as follows.

\vspace{1mm}
\vspace{1mm}
\vspace{1mm}
\vspace{1mm}

\noindent{\bf Theorem 1.2.} {\it
Let $\alpha\in(-1,+\infty)\setminus\mathbb{N}$.
If {\upshape(A)} of Theorem  {\upshape 1.1} holds,  then for any integer $m\geq1$ and any
$s$ large enough,  there exists a  family of solutions $\upsilon_s$
of  problem {\upshape(\ref{1.1})} with $m$ distinct bubbles  accumulating to $q$. Moreover,
$$
\aligned
\lim_{s\rightarrow+\infty}\int_{\Omega}a(x)e^{\upsilon_s}=8\pi(m+1+\alpha)a(q)\phi_1(q).
\endaligned
$$
}

%
%

%
%
%
%
%
%
%

According to Theorems 1.1 and 1.2,  it follows  that
under the assumption (A),
then for any  integer $m\geq1$
there exists a family of solutions  of problem (\ref{1.4}) which exhibits
the phenomenon of $m+1$-bubbling at $q$, namely
$|x-q|^{2\alpha}k(x)e^{-t\phi_1}e^{u_t}\rightharpoonup 8\pi(m+1+\alpha)\delta_{q}$
and $u_t=(m+1+\alpha)G(x,q)+o(1)$. While
for the case $m=0$, by arguing simply along
the initial sketch of the proof of Theorem 1.1,  we  easily
find that problem (\ref{1.4}) always admits a family of solutions
with only one bubble located exactly at singular source  $q$  whether
$q$ is a local maximum point of the functions $\phi_1$ and $a(x)\phi_1$ in
the domain or not.

The proof of  our results  relies on a very
well known  Lyapunov-Schmidt reduction procedure.
The same strategy
has been applied in \cite{DHZ} to
build solutions for the two-dimensional  elliptic problem with
 Hardy-H\'{e}non weight
\begin{equation}\label{1.8}
\aligned
\left\{\aligned
&-\Delta u=|x-q|^{2\alpha}
k(x)e^{-t\phi_1}e^u\,\,\,\,\,
\textrm{in}\,\,\,\,\,\Omega,\\[2mm]
&u=0\,\,\qquad\,\,\qquad\,\,
\quad\quad\qquad\qquad\,\textrm{on}
\,\,\,\partial\Omega,
\endaligned\right.
\endaligned
\end{equation}
as $t\rightarrow+\infty$, where
$\Omega$ is a bounded  smooth  domain in $\mathbb{R}^2$,
$\alpha\in(-1,+\infty)\setminus\mathbb{N}$,
$q\in \Omega $,
$k(x)$ is a given positive smooth function and
$\phi_1$ is a positive eigenfunction of
$-\Delta$ with Dirichlet boundary condition corresponding to the
first eigenvalue $\lambda_1$.
Just like that in equation (\ref{1.4}),
the presence of Hardy-H\'{e}non weight
has significant influence  not only on
the existence of the solution of problem (\ref{1.8})
with a unique bubble at each singular source $q\in\Omega$,
but also on the existence of the solution of  problem (\ref{1.8})
with arbitrarily many   bubbles
 accumulating to some singular source $q\in\Omega$ if
 $q$ is an isolated local maximum point of $\phi_1$.
However, due to the occurrence of Hardy-H\'{e}non weight,
it is necessary to point out that although the anisotropic planar
equation (\ref{1.4}) is seemingly similar to problem (\ref{1.8}),
equation (\ref{1.4})  can  not be  viewed as a special case
of (\ref{1.8}) in higher dimension
even if  the domain has some axial symmetries.
This seems to
imply  that unlike those for solutions of problem (\ref{1.8}) in
\cite{DHZ} whose
multiple clustering bubbles are purely determined by an
isolated local maximum point $q$ of $\phi_1$,
the location  of multiple clustering bubbles in solutions of equation (\ref{1.4}) may be
 not only characterized as isolated local maximum points $q$ of $\phi_1$, but also those of $a(x)\phi_1$,
which needs us to investigate deeply the effect of the interaction between anisotropic coefficient
$a(x)$ and first positive eigenfunction $\phi_1$ on the existence of  solutions with
arbitrarily many bubbles simultaneously  accumulating to singular source
$q$.
This is the delicate description during we carry out the whole  reduction procedure to construct  solutions
of equation (\ref{1.4}) with multiple clustering bubbles around singular source.

\vspace{1mm}
\vspace{1mm}
\vspace{1mm}
\vspace{1mm}

\section{Approximating solutions}
For notational convenience we always fix singular source
$q\in\Omega$  as an isolated local maximum
point of $a(x)\phi_1$
and also a local maximum point of $\phi_1$, and further  assume
\begin{equation}\label{2.1}
\aligned
2\inf_{x\in\overline{B}_{d}(q)}\phi_1(x)>\sup_{x\in \overline{B}_{d}(q)}\phi_1(x)=\phi_1(q)=1,
\endaligned
\end{equation}
where $d>0$ is a small but fixed number, independent of $t$.
For points $\xi=(\xi_1,\ldots,\xi_m)$ with $\xi_i\in \overline{B}_{d}(q)$, we define
\begin{eqnarray}\label{2.2}
\mathcal{O}_t(q):=\left\{\,\xi=(\xi_1,\ldots,\xi_m)\in\big(\overline{B}_d(q)\big)^m\left|\,
a(q)\phi_1(q)-a(\xi_i)\phi_1(\xi_i)\leq\frac{1}{\sqrt{t}},
\,\,\,\,\,
|\xi_i-q|\geq\frac1{t^\beta},
\,\,\,\,\,
|\xi_i-\xi_j|\geq\frac1{t^\beta},
\right.\right.
&&\nonumber\\[0.5mm]
i,j=1,\ldots,m,\,\,\,i\neq j
\Big\},
\qquad\qquad\qquad\qquad\qquad\qquad\qquad\qquad\,
\qquad\qquad\qquad\qquad\qquad\qquad\qquad\,
&&
\end{eqnarray}
where
$\beta$ is given by
\begin{equation}\label{2.3}
\aligned
\beta=\frac{(m+1+\alpha)^2a_2}{2a_1}.
\endaligned
\end{equation}
We thus fix $\xi\in\mathcal{O}_t(q)$. For numbers $\mu_0>0$ and
$\mu_i>0$, $i=1,\ldots,m$, yet to
be chosen, we define
\begin{equation}\label{2.4}
\aligned
u_0(x)=\log
\frac{8\mu_0^2(1+\alpha)^2}{k(q)(\varepsilon_{0}^2\mu_{0}^2+|x-q|^{2(1+\alpha)})^2},
\qquad\qquad
u_i(x)=\log
\frac{8\mu_i^2}{k(\xi_i)|\xi_i-q|^{2\alpha}(\varepsilon_{i}^2\mu_{i}^2+|x-\xi_i|^2)^2},
\endaligned
\end{equation}
which, respectively, solve
\begin{equation}\label{2.5}
\aligned
-\Delta
u_0=\varepsilon_0^2k(q)|x-q|^{2\alpha}e^{u_0}
\qquad
\textrm{in}
\quad
\mathbb{R}^2,
\,\qquad\quad\quad\,
\int_{\mathbb{R}^2}
\varepsilon_0^2k(q)|x-q|^{2\alpha}e^{u_0}=8\pi(1+\alpha),
\endaligned
\end{equation}
and
\begin{equation}\label{2.6}
\aligned
-\Delta
u_i=\varepsilon_i^2k(\xi_i)|\xi_i-q|^{2\alpha}e^{u_i}
\quad\,\,
\textrm{in}
\quad
\mathbb{R}^2,
\,\qquad\quad\quad\,
\int_{\mathbb{R}^2}\varepsilon_i^2k(\xi_i)|\xi_i-q|^{2\alpha}e^{u_i}=8\pi,
\endaligned
\end{equation}
where
\begin{equation}\label{2.7}
\aligned
\varepsilon_0=\varepsilon_0(t)\equiv e^{-\frac12t},
\,\qquad\qquad\,
\varepsilon_i=\varepsilon_i(t)\equiv e^{-\frac12t\phi_1(\xi_i)}.
\endaligned
\end{equation}
Furthermore, we set
\begin{equation}\label{2.7a}
\aligned
\rho_0:=\varepsilon_0^{\frac{1}{1+\alpha}}=\exp\left\{
-\frac1{2(1+\alpha)}t\right\},
\,\,\,\quad\,\,\,\,
v_0:=\mu_0^{\frac{1}{1+\alpha}},
\,\,\,\quad\,\,\,
\gamma_i:=\frac{1}{\varepsilon_0}\varepsilon_i\mu_i=\mu_i\exp\left\{
-\frac12t\big[\phi_1(\xi_i)-1\big]\right\}.
\endaligned
\end{equation}

We defined the approximate solution  of problem (\ref{1.4}) by
\begin{equation}\label{2.8}
\aligned
U(x):=\sum\limits_{i=0}^{m}\,U_i(x)
=\sum\limits_{i=0}^{m}\,
\big[u_i(x)+H_i(x)\big],
\endaligned
\end{equation}
where $H_i(x)$ is a correction term defined as the solution of
\begin{equation}\label{2.9}
\left\{\aligned
&\Delta_{a}
H_i+\nabla\log a(x)\nabla u_i=0
\,\,\,\,
\textrm {in}\,\,\,\,\,\Omega,\\[2mm]
&H_i=-u_i
\,\quad\,\,\,\,\,\quad\quad
\quad\quad\quad\quad
\textrm{on}\,\,\,\po.
\endaligned\right.
\end{equation}

\vspace{1mm}

\noindent{\bf Lemma 2.1.}\,\,{\it For
any $0<\sigma<\min\{1,\,2(1+\alpha)\}$, any points
$\xi=(\xi_1,\ldots,\xi_m)\in\mathcal{O}_t(q)$
and any $t$ large enough,
\begin{equation}\label{2.10}
\aligned
H_0(x)=(1+\alpha)H(x,q)-\log
\frac{8\mu_0^2(1+\alpha)^2}{k(q)}+O\left(\rho_0^\sigma v_0^\sigma\right),
\endaligned
\end{equation}
\begin{equation}\label{2.11}
\aligned
H_i(x)=H(x,\xi_i)-\log
\frac{8\mu_i^2}{k(\xi_i)|\xi_i-q|^{2\alpha}}+O\left(\varepsilon_i^\sigma\mu_i^\sigma\right),
\,\quad\,i=1,\ldots,m,
\endaligned
\end{equation}
uniformly in $\overline{\Omega}$, where $H$ is the regular part of Green's function defined in
{\upshape(\ref{1.3})}.
}

\vspace{1mm}

\begin{proof}
According to Lemma 2.4 in \cite{WYZ}, we need only to prove the expansion for  $H_0$.
Observe that, on  $\po$,
$$
\aligned
H_0(x)=-u_0(x)=4(1+\alpha)\log|x-q|-\log\frac{8\mu_0^2(1+\alpha)^2}{k(q)}+O\left(\varepsilon_0^2\mu_0^2\right).
\endaligned
$$
From (\ref{1.2})-(\ref{1.3}) it follows that the regular part of Green's function $H(x,y)$ satisfies
$$
\aligned
\left\{\aligned
&-\Delta_a H(x,y)=-4\frac{(x-y)\cdot\nabla\log a(x)}{|x-y|^2},
\,\,\,\,\,\,
x\in\Omega,\\
&H(x,y)=4\log|x-y|,
\,\qquad\qquad\qquad\,\,\,\,\quad\,
x\in\partial\Omega.
\endaligned\right.
\endaligned
$$
If we set $Z_0(x)=H_0(x)-(1+\alpha)H(x,q)+\log\frac{8\mu_0^2(1+\alpha)^2}{k(q)}$,
then $Z_0(x)$ satisfies
$$
\aligned
\left\{\aligned
&-\Delta_a Z_0(x)=
\frac{4(1+\alpha)\varepsilon_{0}^2\mu_{0}^2}{\varepsilon_{0}^2\mu_{0}^2+|x-q|^{2(1+\alpha)}}
\frac{(x-q)\cdot\nabla\log a(x)}{|x-q|^2}
\,\ \ \,\,\textrm{in}\,\,\,\,\,\Omega,\\[2mm]
&Z_0(x)=O\left(\varepsilon_0^2\mu_0^2\right)\,\
\quad\qquad\qquad\qquad\qquad\quad\,
\qquad\qquad\qquad\,\textrm{on}\,\,\,\po.
\endaligned\right.
\endaligned
$$
Using the polar coordinates with center $q$, i.e. $r=|x-q|$ and
making  the change of  variables $s=r/(\rho_0v_0)$, we
have that for any  $1<p<2$,
$$
\aligned
\int_{\Omega}\left|
\frac{4(1+\alpha)\varepsilon_{0}^2\mu_{0}^2}{\varepsilon_{0}^2\mu_{0}^2+|x-q|^{2(1+\alpha)}}
\frac{(x-q)\cdot\nabla\log a(x)}{|x-q|^2}
\right|^p dx
&\leq C\int_{0}^{\diam(\Omega)}\left|\frac{\varepsilon_{0}^2\mu_{0}^2}{r\big(\varepsilon_{0}^2\mu_{0}^2+r^{2(1+\alpha)}\big)}
\right|^p rdr\\
&\leq C(\rho_0v_0)^{2-p}\int_{0}^{\diam(\Omega)/(\rho_0v_0)}
\frac{s^{1-p}}{\big(1+s^{2(1+\alpha)}\big)^p}
ds\\[1mm]
&\leq C\big[(\rho_0v_0)^{2-p}
+(\rho_0v_0)^{2p(1+\alpha)}
\big].
\endaligned
$$
Applying $L^p$ theory,
$$
\aligned
\left\|Z_0\right\|_{W^{2,p}(\Omega)}\leq
C\big(\big\|-\Delta_a
Z_0\big\|_{L^p(\Omega)}+\left\|Z_0\right\|_{C^{2}(\partial\Omega)}\big)
\leq C\big[(\rho_0v_0)^{(2-p)/p}
+(\rho_0v_0)^{2(1+\alpha)}
\big].
\endaligned
$$
By Sobolev embedding, we obtain that for any $0<\gamma<2-(2/p)$,
$$
\aligned
\left\|Z_0\right\|_{C^\gamma(\overline{\Omega})}\leq C\big[(\rho_0v_0)^{(2-p)/p}
+(\rho_0v_0)^{2(1+\alpha)}
\big].
\endaligned
$$
This proves expansion (\ref{2.10}) with $\sigma<\min\{(2-p)/p,\,2(1+\alpha)\}$.
\end{proof}

\vspace{1mm}
\vspace{1mm}
\vspace{1mm}

Let us write
$$
\aligned
\Omega_t=\varepsilon_0^{-1}\Omega=e^{\frac12t}\Omega,
\qquad\quad
q'=q/\varepsilon_0
\qquad\quad
\textrm{and}
\qquad\quad
\xi'_i=\xi_i/\varepsilon_0,
\quad
i=1,\ldots,m.
\endaligned
$$
Then $u(x)$ solves  problem  (\ref{1.4}) if and only if $\omega(y)\equiv u(\varepsilon_0 y)-2t$ satisfies
\begin{equation}\label{2.17}
\aligned
\left\{\aligned
&-\Delta_{a(\varepsilon_0 y)}
\omega=|\varepsilon_0y-q|^{2\alpha}\kappa(y,t)e^\omega
\quad
\textrm{in}\,\,\,\,\,\Omega_t,\\[2mm]
&\omega=-2t
\qquad\qquad\qquad\qquad\,\,\,
\qquad\quad\,\,\,\textrm{on}
\,\,\,\partial\Omega_t,
\endaligned\right.
\endaligned
\end{equation}
where
\begin{equation}\label{2.18}
\aligned
\kappa(y,t)\equiv k(\varepsilon_0 y)\exp
\left\{-t\big[\phi_1(\varepsilon_0 y)-1\big]\right\}.
\endaligned
\end{equation}
Let us  define the initial approximate solution of (\ref{2.17}) as
\begin{equation}\label{2.19}
\aligned
V(y)=U(\varepsilon_0 y)-2t
\endaligned
\end{equation}
with $U$  given by (\ref{2.8}). Hence if
we try to look for a
solution  of equation (\ref{2.17}) in the form
$\omega=V+\phi$ with $\phi$ a lower order correction, then
  (\ref{2.17}) can be stated as to find $\phi$ a solution
of
\begin{equation}\label{2.34}
\left\{\aligned
&\mathcal{L}(\phi):=-\Delta_{a(\varepsilon_0 y)} \phi-W\phi=E+N(\phi)
\,\,\,\,
\textrm{in}\,\,\,\,\,\Omega_t,\\[2mm]
&\phi=0
\quad\qquad\qquad\qquad\,\
\qquad\qquad\qquad\quad\,\,
\textrm{on}\,\,\,
\partial\Omega_t,
\endaligned\right.
\end{equation}
where
\begin{equation}\label{2.20}
\aligned
W(y)=|\varepsilon_0y-q|^{2\alpha}\kappa(y,t)e^{V},
\endaligned
\end{equation}
 the ``error term'' is
\begin{equation}\label{2.21}
\aligned
E(y)=\Delta_{a(\varepsilon_0 y)} V+|\varepsilon_0y-q|^{2\alpha}\kappa(y,t)e^{V}.
\endaligned
\end{equation}
and the ``nonlinear term''  is given by
\begin{equation}\label{2.35}
\aligned
N(\phi)=|\varepsilon_0y-q|^{2\alpha}\kappa(y,t)e^{V}\big(e^\phi-1-\phi\big),
\endaligned
\end{equation}
In order to understand how well $V(y)$ solves equation (\ref{2.17}) so that the
``error term''
$E(y)$ is sufficiently small near $q$ and each $\xi_i$ with $i=1,\ldots,m$, a delicate
ingredient is to make the following
precise choices of  the concentration parameters $\mu_0$ and $\mu_i$:
\begin{equation}\label{2.12}
\aligned
\log
\frac{8\mu_0^2(1+\alpha)^2}{k(q)}=
(1+\alpha)H(q,q)+\sum_{j=1}^m G(q,\xi_j),
\endaligned
\end{equation}
\begin{equation}\label{2.13}
\aligned
\log
\frac{8\mu_i^2}{k(\xi_i)|\xi_i-q|^{2\alpha}}=
H(\xi_i,\xi_i)+(1+\alpha)G(\xi_i,q)
+\sum_{j=1,\,j\neq
i}^m G(\xi_i,\xi_j),
\,\quad\,i=1,\ldots,m.
\endaligned
\end{equation}
Here, $\mu_0$ and $\mu_i$  are {\it a priori}  functions of $\xi$ in $\mathcal{O}_t(q)$ and
hence $\mu_0=\mu_0(\xi)$ and $\mu_i=\mu_i(\xi)$ for all $i=1,\ldots,m$.
Thanks to the definition of $\mathcal{O}_t(q)$ in (\ref{2.2}),
there exists a constant $C>0$ independent of $t$ such that
\begin{equation}\label{2.14}
\aligned
\frac1C\leq\mu_0\leq Ct^{2m\beta}
\quad\ \quad\,\,\textrm{and}\quad\quad\,\,
\big|\partial_{\xi_{kl}}\log\mu_0\big|\leq Ct^{\beta},\,
\,\quad\ \,\forall\,\,k=1,\ldots,m,\,\,l=1,2,
\endaligned
\end{equation}
and
\begin{equation}\label{2.15}
\aligned
\frac1C\leq\mu_i\leq Ct^{(2m+\alpha)\beta}
\quad\ \quad\,\,\textrm{and}\quad\quad\,\,
\big|\partial_{\xi_{kl}}\log\mu_i\big|\leq Ct^{\beta},\,
\,\quad\,\forall\,\,i, k=1,\ldots,m,\,\,l=1,2.
\endaligned
\end{equation}

Finally, we claim that the following behavior for $E(y)$ holds:
for any $0<\sigma<\min\{1/2,\,1-1/(2\beta),\,2(1+\alpha)\}$ it yields that
if   $|y-q'|\leq1/(\varepsilon_0 t^{2\beta})$,
\begin{eqnarray}\label{2.31}
E(y)=\left(\frac{\varepsilon_0}{\rho_0v_0}\right)^2\frac{8(1+\alpha)^2\big|\frac{\varepsilon_0 y-q}{\rho_0v_0}\big|^{2\alpha}}{\big(1+\big|\frac{\varepsilon_0 y-q}{\rho_0v_0}\big|^{2(1+\alpha)}\big)^2}
O\left(\varepsilon_0^\sigma|y-q'|^{\sigma}
+\rho_0^\sigma v_0^\sigma
+\sum_{j=1}^m
\varepsilon_j^\sigma\mu_j^\sigma
\right)
+\sum\limits_{j=1}^{m}O\left(
\varepsilon_0^2\varepsilon_j^2\mu_j^2t^{4\beta}
\right),
\end{eqnarray}
and if  $|y-\xi'_i|\leq1/(\varepsilon_0 t^{2\beta})$
with some $i\in\{1,\ldots,m\}$,
\begin{eqnarray}\label{2.29}
E(y)=\frac{1}{\gamma_{i}^2}\frac{8}{\big(1+\big|\frac{y-\xi'_i}{\gamma_{i}}\big|^2\big)^2}
O\left(\varepsilon_0^\sigma|y-\xi_i'|^{\sigma}
+
\rho_0^\sigma v_0^\sigma
+\sum_{j=1}^m
\varepsilon_j^\sigma\mu_j^\sigma\right)
+O\left(
\varepsilon_0^4\mu_0^2t^{(4+2\alpha)\beta}
+\sum\limits_{j=1,\,j\neq i}^{m}
\varepsilon_0^2\varepsilon_j^2\mu_j^2t^{4\beta}
\right),
\end{eqnarray}
while if  $|y-q'|>1/(\varepsilon_0 t^{2\beta})$ and
$|y-\xi'_i|>1/(\varepsilon_0 t^{2\beta})$ for all
 $i=1,\ldots,m$,
\begin{equation}\label{2.33}
\aligned
E(y)=O\left(\frac{\varepsilon_0^2e^{-t\phi_1(\varepsilon_0 y)}}{\,|\varepsilon_0 y-q|^{4+2\alpha}\,}\prod_{i=1}^m\frac{1}{\,|\varepsilon_0 y-\xi_i|^4\,}
\right)+
O\left(
\varepsilon_0^4\mu_0^2t^{(8+4\alpha)\beta}
\right)
+\sum\limits_{i=1}^{m}O\left(
\varepsilon_0^2\varepsilon_i^2\mu_i^2t^{8\beta}
\right).
\endaligned
\end{equation}
In fact, recalling   that $E(y)=\Delta_{a(\varepsilon_0 y)} V+W$ with $V$ and $W$ given by (\ref{2.19}) and (\ref{2.20}), respectively,
we first have
\begin{eqnarray}\label{2.22}
-\Delta_{a(\varepsilon_0 y)} V(y)
=\varepsilon_0^4k(q)|x-q|^{2\alpha}e^{u_0}
+\varepsilon_0^2\sum\limits_{i=1}^{m}\varepsilon_i^2k(\xi_i)|\xi_i-q|^{2\alpha}e^{u_i}
\,\qquad\qquad\qquad\,\,&&\nonumber\\
=\frac{8\varepsilon_0^4\mu_0^2(1+\alpha)^2|x-q|^{2\alpha}}{\,(\varepsilon_{0}^2\mu_{0}^2+|x-q|^{2(1+\alpha)})^2\,}
+\sum\limits_{i=1}^{m}
\frac{8\varepsilon_0^2\varepsilon_i^2\mu_i^2}{\,(\varepsilon_{i}^2\mu_{i}^2+|x-\xi_i|^2)^2\,}
\,\qquad\qquad\,\,\,\,
&&\nonumber\\
=\left(\frac{\varepsilon_0}{\rho_0v_0}\right)^2\frac{8(1+\alpha)^2\big|\frac{\varepsilon_0 y-q}{\rho_0v_0}\big|^{2\alpha}}{\,\big(1+\big|\frac{\varepsilon_0 y-q}{\rho_0v_0}\big|^{2(1+\alpha)}\big)^2\,}+\sum\limits_{i=1}^{m}
\frac{1}{\gamma_{i}^2}\frac{8}{\big(1+\big|\frac{y-\xi'_i}{\gamma_{i}}\big|^2\big)^2}.
\,\quad\quad\,
&&
\end{eqnarray}
For the expression  $W$ near singular source $q$, we can compute
\begin{eqnarray*}\label{2.27}
W(y)=\varepsilon_0^4|\varepsilon_0y-q|^{2\alpha}\kappa(y,t)\exp
\left\{
\sum\limits_{j=0}^{m}\,
\big[u_j(\varepsilon_0 y)+H_j(\varepsilon_0 y)\big]\right\}
\quad\qquad\qquad
\qquad\,\,\,
\quad\quad\,\,\,\,
&&\nonumber\\
=\varepsilon_0^4k(\varepsilon_0 y)|\varepsilon_0y-q|^{2\alpha}
\exp\left\{-t\big[\phi_1(\varepsilon_0 y)-1\big]
+\log
\frac{8\mu_0^2(1+\alpha)^2}{k(q)(\varepsilon_{0}^2\mu_{0}^2+|x-q|^{2(1+\alpha)})^2}
\right.
&&\nonumber\\
\left.
+H_0(\varepsilon_0 y)
+
\sum\limits_{j=1}^{m}\,
\big[u_j(\varepsilon_0 y)
+H_j(\varepsilon_0 y)\big]\right\}
\qquad\qquad\qquad\qquad\quad
\qquad\qquad\qquad\,\,\,\,\,\,
&&\nonumber\\
=\left(\frac{\varepsilon_0}{\rho_0v_0}\right)^2\frac{8(1+\alpha)^2\big|\frac{\varepsilon_0 y-q}{\rho_0v_0}\big|^{2\alpha}}{\,\big(1+\big|\frac{\varepsilon_0 y-q}{\rho_0v_0}\big|^{2(1+\alpha)}\big)^2\,}
\times\frac{k(\varepsilon_0 y)}{k(q)}
\times\exp
\Big\{-t\big[\phi_1(\varepsilon_0 y)-\phi_1(q)\big]
\Big\}
\quad\,\,
&&\nonumber\\
\times\exp\left\{H_0(\varepsilon_0 y)
+
\sum\limits_{j=1}^{m}\,
\big[u_j(\varepsilon_0 y)
+H_j(\varepsilon_0 y)\big]\right\}.
\qquad\qquad\quad\,\,\,\,
\qquad\qquad\qquad\quad\,
&&
\end{eqnarray*}
Using (\ref{2.4}), (\ref{2.10}), (\ref{2.11}) and the fact that $H(\cdot,x)$ is $C^\sigma(\Omega)$
for any $x\in\Omega$ and $0<\sigma<\min\{1/2,\,1-1/(2\beta),\,2(1+\alpha)\}$,
we obtain that  for $|y-q'|\leq1/(\varepsilon_0 t^{2\beta})$,
$$
\aligned
H_0&(\varepsilon_0 y)
+
\sum\limits_{j=1}^{m}\,
\big[u_j(\varepsilon_0 y)
+H_j(\varepsilon_0 y)\big]\\
&=(1+\alpha)H(\varepsilon_0 y,q)-\log
\frac{8\mu_0^2(1+\alpha)^2}{k(q)}+O\left(\rho_0^\sigma v_0^\sigma\right)\\
&\quad\,+
\sum\limits_{j=1}^{m}\,
\left[\log\frac{1}{|q-\xi_j|^4}+H(\varepsilon_0 y,\xi_j)
+O\left(
\frac{|\varepsilon_0 y-q|^2+
2\langle \varepsilon_0 y-q, q-\xi_j\rangle
+\varepsilon_j^2\mu_j^2}{|q-\xi_j|^2}
\right)+O\left(\varepsilon_j^\sigma\mu_j^\sigma\right)
\right]\\
&=(1+\alpha)H(q,q)
-\log\frac{8\mu_0^2(1+\alpha)^2}{k(q)}
+\sum_{j=1}^m G(q,\xi_j)
+
O\left(\varepsilon_0^\sigma|y-q'|^{\sigma}\right)
+O\left(\rho_0^\sigma v_0^\sigma\right)+\sum_{j=1}^mO\left(\varepsilon_j^\sigma\mu_j^\sigma\right)\\
&=O\left(\varepsilon_0^\sigma|y-q'|^{\sigma}\right)
+O\left(\rho_0^\sigma v_0^\sigma\right)+\sum_{j=1}^mO\left(\varepsilon_j^\sigma\mu_j^\sigma\right),
\endaligned
$$
where the last equality is due to the choice of $\mu_0$
in (\ref{2.12}).
Therefore if   $|y-q'|\leq1/(\varepsilon_0 t^{2\beta})$,
\begin{equation}\label{2.30}
\aligned
W(y)
=\left(\frac{\varepsilon_0}{\rho_0v_0}\right)^2\frac{8(1+\alpha)^2\big|\frac{\varepsilon_0 y-q}{\rho_0v_0}\big|^{2\alpha}}
{\big(1+\big|\frac{\varepsilon_0 y-q}{\rho_0v_0}\big|^{2(1+\alpha)}\big)^2}
\left[1
+O\left(\varepsilon_0^\sigma|y-q'|^{\sigma}\right)
+O\left(\rho_0^\sigma v_0^\sigma\right)+\sum_{j=1}^mO\left(\varepsilon_j^\sigma\mu_j^\sigma\right)
\right].
\endaligned
\end{equation}
Similarly, if
$|y-\xi'_i|\leq1/(\varepsilon_0 t^{2\beta})$ for some $i\in\{1,\ldots,m\}$,
\begin{eqnarray}\label{2.28}
W(y)=\frac{1}{\gamma_{i}^2}\frac{8}{\big(1+\big|\frac{y-\xi'_i}{\gamma_{i}}\big|^2\big)^2}
\left[1+
O\left(\varepsilon_0^\sigma|y-\xi'_i|^{\sigma}\right)
+O\left(\rho_0^\sigma v_0^\sigma\right)+\sum_{j=1}^mO\left(\varepsilon_j^\sigma\mu_j^\sigma\right)\right],
\end{eqnarray}
while if  $|y-q'|>1/(\varepsilon_0 t^{2\beta})$ and
$|y-\xi'_i|>1/(\varepsilon_0 t^{2\beta})$
for all $i=1,\ldots,m$,
\begin{equation}\label{2.32}
\aligned
W(y)
=O\left(\frac{\varepsilon_0^2e^{-t\phi_1(\varepsilon_0 y)}}
{\,|\varepsilon_0 y-q|^{4+2\alpha}\,}\prod_{i=1}^m\frac{1}{\,|\varepsilon_0 y-\xi_i|^4\,}
\right).
\endaligned
\end{equation}
These combined with (\ref{2.22})  imply that the expansions of $E(y)$ in (\ref{2.31})-(\ref{2.33})  hold.

\vspace{1mm}
\vspace{1mm}
\vspace{1mm}

\section{The linearized problem}
In this section we solve
the following linear problem: given $h\in L^{\infty}(\Omega_t)$
and points $\xi=(\xi_1,\ldots,\xi_m)\in\mathcal{O}_t(q)$,
we find a function $\phi$ such that for certain  scalars $c_{ij}$,
$i=1,\ldots,m$, $j=1,2$, one has
\begin{equation}\label{3.1}
\aligned
\left\{
\aligned
&\mathcal{L}(\phi)=-\Delta_{a(\varepsilon_0 y)} \phi-W\phi=h+\frac1{a(\varepsilon_0 y)}\sum\limits_{i=1}^m\sum\limits_{j=1}^2c_{ij}\chi_iZ_{ij}
\quad\,\,\,\,
\textrm{in}\,\,\,\,\,\Omega_t,\\
&\phi=0
\quad\quad\quad\quad\quad
\,\quad\quad\qquad\quad\,\,\,\,\,
\qquad\qquad\qquad\quad\qquad\ \,
\quad\quad\,\,
\textrm{on}
\,\,\,\po_t,\\[1mm]
&\int_{\Omega_t}\chi_iZ_{ij}\phi=0
\qquad\qquad\qquad\qquad\qquad\quad\ \,
\forall\,\,i=1,\ldots,m,
\,\,j=1,2,
\endaligned
\right.
\endaligned
\end{equation}
where $W=|\varepsilon_0y-q|^{2\alpha}\kappa(y,t)e^{V}$
satisfies (\ref{2.30})-(\ref{2.32}), and $Z_{ij}$, $\chi_i$ are defined as follows:
let $R_0$ be  a large but fixed positive number  and $\chi(r)$ be
a radial  smooth non-increasing  cut-off function  satisfying
$0\leq\chi(r)\leq1$, $\chi(r)=1$ for $r\leq R_0$ and $\chi(r)=0$ for $r\geq R_0+1$.
Set
\begin{equation}\label{3.2}
\aligned
\mathcal{Z}_{q}(z)=\frac{|z|^{2(1+\alpha)}-1}{|z|^{2(1+\alpha)}+1},
\,\,\quad\quad\quad\,\,
\mathcal{Z}_{0}(z)=\frac{|z|^2-1}{|z|^2+1},
\,\,\quad\quad\quad\,\,
\mathcal{Z}_{j}(z)=\frac{4z_j}{|z|^2+1},\,\,\,\,j=1,\,2.
\endaligned
\end{equation}
Then we define
\begin{equation}\label{3.5}
\aligned
\chi_q(y)=\chi\left(\frac{\big|\varepsilon_0y-q\big|}{\rho_0 v_0}\right)
\,\,\qquad\,\,\ \,\textrm{and}
\,\,\qquad\,\,\,
Z_{q}(y)=\frac{\varepsilon_0}{\rho_0 v_0}\mathcal{Z}_q\left(\frac{\varepsilon_0y-q}{\rho_0 v_0}\right),
\endaligned
\end{equation}
and for any $i=1,\ldots,m$ and $j=0,1,2$,
\begin{equation}\label{3.6}
\aligned
\chi_i(y)=\chi\left(\frac{\big|y-\xi_i'\big|}{\gamma_i}\right)
\,\,\,\qquad\,\,\,
\textrm{and}
\,\,\,\qquad\,\,\,
Z_{ij}(y)=\frac{1}{\gamma_i}\mathcal{Z}_j\left(\frac{y-\xi_i'}{\gamma_i}\right).
\endaligned
\end{equation}

Equation (\ref{3.1}) will be solved for
$h\in L^\infty(\Omega_t)$, but we need to estimate the size of
the solution by introducing the following  norm:
\begin{eqnarray}\label{3.7}
\quad
\|h\|_{*}:=
\left\|\left[
\varepsilon_0^2+
\left(\frac{\varepsilon_0}{\rho_0v_0}\right)^2\frac{\big|\frac{\varepsilon_0 y-q}{\rho_0v_0}\big|^{2\alpha}}
{\,\big(1+\big|\frac{\varepsilon_0 y-q}{\rho_0v_0}\big|\big)^{4+2\hat{\alpha}+2\alpha}\,}
+
\sum\limits_{i=1}^m\frac{1}{\gamma_i^2}\frac{1}{\big(1+\big|\frac{y-\xi'_i}{\gamma_i}\big|\big)^{4+2\hat{\alpha}}}
\right]^{-1}h(y)
\right\|_{L^\infty(\Omega_t)},
\end{eqnarray}
where $\hat{\alpha}+1$ is a  small but fixed positive number, independent of $t$,
such that $-1<\hat{\alpha}<\min\big\{\alpha,\,-2/3\big\}$.

\vspace{1mm}
\vspace{1mm}
\vspace{1mm}
\vspace{1mm}

\noindent{\bf Proposition 3.1.}\,\,{\it
Let $m$ be a positive  integer. Then
there exist constants $t_m>1$ and $C>0$ such
that for any $t>t_m$,
any points $\xi=(\xi_1,\ldots,\xi_m)\in\mathcal{O}_t(q)$ and any
$h\in L^{\infty}(\Omega_t)$, there is a unique solution
$\phi=\mathcal{T}(h)$ and scalars
$c_{ij}$,
$i=1,\ldots,m$, $j=1,2$
to problem {\upshape(\ref{3.1})}, which satisfies
\begin{equation}\label{3.8}
\aligned
\|\mathcal{T}(h)\|_{L^{\infty}(\Omega_t)}\leq Ct\|h\|_{*}
\qquad\qquad
\textrm{and}
\qquad\qquad
|c_{ij}|\leq C\gamma_i^{-1}\|h\|_{*}.
\endaligned
\end{equation}

}
\begin{proof}
The proof will be divided into  four steps which we state and prove next.

\vspace{1mm}

{\bf Step 1:}  Building a suitable  barrier defined in
$\widetilde{\Omega}_t:=
\Omega_t\setminus\big[\bigcup_{i=1}^mB_{R_1\gamma_i}(\xi'_i)
\cup B_{R_1\rho_0 v_0/\varepsilon_0}(q')\cup B^c_{2d/\varepsilon_0}(q')\big]$
for some $R_1$ large but $d$ small, independent of $t$.

\vspace{1mm}
\vspace{1mm}
\vspace{1mm}

\noindent{\bf Lemma 3.1.}\,\,{\it There exist positive constants $R_1$ and $C$,
independent of $t$, such that
for  any points $\xi=(\xi_1,\ldots,\xi_m)\in\mathcal{O}_t(q)$ and any
$t$ large enough, there exists
$\psi:\overline{\widetilde{\Omega}}_t\rightarrow\mathbb{R}$
smooth  and
positive verifying
$$
\aligned
\mathcal{L}(\psi)=-\Delta_{a(\varepsilon_0 y)} \psi-W\psi\geq
\left(\frac{\varepsilon_0}{\rho_0v_0}\right)^2\frac{1}
{\big|\frac{\varepsilon_0 y-q}{\rho_0v_0}\big|^{4+2\hat{\alpha}}}
+
\sum_{i=1}^m\frac{1}{\gamma_i^2}\frac{1}{\big|\frac{y-\xi'_i}{\gamma_i}\big|^{4+2\hat{\alpha}}}
+\varepsilon_0^2
\,\,\,\quad\,\,
\textrm{in}\,\,\,\,\widetilde{\Omega}_t.
\endaligned
$$
Moreover, $\psi$ is uniformly  bounded, i.e.
$$
\aligned
1<\psi\leq C\,\,\quad\textrm{in}\,\,\,
\,\overline{\widetilde{\Omega}}_t.
\endaligned
$$
}

\begin{proof}
Let us take
$$
\aligned
\psi=C_1\left(\Psi_0(\varepsilon_0 y)-\frac{1}
{\,\big|\frac{\varepsilon_0 y-q}{\rho_0v_0}\big|^{2(1+\hat{\alpha})}\,}\right)
+C_1\sum_{i=1}^m
\left(\Psi_0(\varepsilon_0 y)-
\frac{1}{\,\big|\frac{y-\xi'_i}{\gamma_i}\big|^{2(1+\hat{\alpha})}\,}
\right),
\endaligned
$$
where $\Psi_0$ satisfies
$-\Delta_a\Psi_0=1$ in $\Omega$, $\Psi_0=2$ on $\partial\Omega$.
Since $\Psi_0\geq2$ in $\Omega$,
it is directly checked that, choosing the positive
constant $C_1$ larger if necessary, $\psi$  meets all the
conditions of the lemma for numbers $R_1$ and $t$ large enough.
\end{proof}

\vspace{1mm}

{\bf Step 2:} An auxiliary linear equation.  Given $h\in L^{\infty}(\Omega_t)$
and $\xi=(\xi_1,\ldots,\xi_m)\in\mathcal{O}_t(q)$, we first study the linear equation
\begin{equation}\label{3.10}
\aligned
\left\{\aligned
&\mathcal{L}(\phi)=-\Delta_{a(\varepsilon_0 y)}\phi-W\phi=h
\,\,\quad\,\textrm{in}\,\,\,\,\,\Omega_t,\\[2mm]
&\phi=0
\,\,\,\,\qquad\qquad\,
\,\,\,\,\quad\,
\,\quad\,\,\,\,
\,\,\,\quad\,\,\,
\quad\,\textrm{on}
\,\,\,\po_t.
\endaligned\right.
\endaligned
\end{equation}
For solutions of (\ref{3.10}) involving more
orthogonality
conditions than those in (\ref{3.1}), we have
the following a priori estimate.

\vspace{1mm}
\vspace{1mm}
\vspace{1mm}
\vspace{1mm}

\noindent{\bf Lemma 3.2.}\,\,{\it There exist $R_0>0$ and  $t_m>1$ such that for any $t>t_m$
and any solution $\phi$ of {\upshape (\ref{3.10})} with the orthogonality conditions
\begin{equation}\label{3.14}
\aligned
\int_{\Omega_t}\chi_qZ_q\phi=0
\,\qquad\,
\textrm{and}
\,\qquad\,\int_{\Omega_t}\chi_iZ_{ij}\phi=0,\,\,\,\,\,
\,\,\,\,\,i=1,\ldots,m,\,\,\,j=0,1,2,
\endaligned
\end{equation}
we have
\begin{equation}\label{3.12}
\aligned
\|\phi\|_{L^{\infty}(\Omega_t)}\leq C
\|h\|_{*},
\endaligned
\end{equation}
where $C>0$ is independent of $t$.}

\vspace{1mm}

\begin{proof}
Take $R_0=2R_1$ with $R_1$  the constant in the previous
 step.
Since $\xi=(\xi_1,\ldots,\xi_m)\in\mathcal{O}_t(q)$,
$\rho_0 v_0=o(1/t^{2\beta})$ and
$\varepsilon_0\gamma_i=o(1/t^{2\beta})$ for $t$ large enough,
we have $B_{R_1\rho_0 v_0/\varepsilon_0}(q')$ and
$B_{R_1\gamma_i}(\xi'_i)$, $i=1,\ldots,m$, disjointed and included in $\Omega_t$.
Recalling  the  barrier $\psi$ in the previous lemma, we first claim that
the operator $\mathcal{L}$ satisfies the maximum principle in
$\widetilde{\Omega}_t$,
namely if $\phi$ is a  supersolution of   $\mathcal{L}(\phi)=-
\Delta_{a(\varepsilon_0 y)}\phi-W\phi\geq0$ in $\widetilde{\Omega}_t$,
$\phi\geq0$ on $\partial\widetilde{\Omega}_t$,
then $\phi\geq0$ in $\widetilde{\Omega}_t$.
In fact, suppose by contradiction that the operator $\mathcal{L}$ does not satisfy the maximum principle in
$\widetilde{\Omega}_t$.  Since $\psi>0$ in $\widetilde{\Omega}_t$,
the function $\phi/\psi$ has a  negative minimum  point $y_0$
in $\widetilde{\Omega}_t$. A simple computation deduces
$$
\aligned
-\Delta_{a(\varepsilon_0 y)}\left(\frac{\phi}{\psi}\right)=\frac{1}{\psi^2}\big[
\psi\mathcal{L}(\phi)
-\phi\mathcal{L}(\psi)
\big]+\frac{2}{\psi}\nabla \psi\nabla\left(\frac{\phi}{\psi}\right).
\endaligned
$$
This, together with the fact that $\mathcal{L}(\psi)>0$ in $\widetilde{\Omega}_t$, gives  $-\Delta_{a(\varepsilon_0 y)}\big(\phi/\psi\big)(y_0)>0$.
 On the other hand,
$$
\aligned
-\Delta\left(\frac{\phi}{\psi}\right)=-\Delta_{a(\varepsilon_0 y)}\left(\frac{\phi}{\psi}\right)
+
\varepsilon_0\nabla \log a(\varepsilon_0y)\nabla\left(\frac{\phi}{\psi}\right).
\endaligned
$$
Then $-\Delta \big(\phi/\psi\big)(y_0)>0$,  which contradicts to the fact that
$y_0$ is a minimum  point of $\phi/\psi$ in $\widetilde{\Omega}_t$.

Let $h$ be bounded and $\phi$ be a  solution to (\ref{3.10}) satisfying (\ref{3.14}).
We define  the ``inner norm'' of $\phi$ as
\begin{equation}\label{3.11}
\aligned
\|\phi\|_{**}=\sup_{y\in
\bigcup_{i=1}^mB_{R_1\gamma_i}(\xi'_i)\cup B_{R_1\rho_0 v_0/\varepsilon_0}(q')
\cup\big(\Omega_t\setminus B_{2d/\varepsilon_0}(q')\big)}
|\phi(y)|,
\endaligned
\end{equation}
and  claim that there is a constant $C>0$ independent of  $t$ such that
for  any points $\xi=(\xi_1,\ldots,\xi_m)\in\mathcal{O}_t(q)$,
\begin{equation}\label{3.12}
\aligned
\|\phi\|_{L^{\infty}(\Omega_t)}\leq C\left(\|\phi\|_{**}+
\|h\|_{*}\right).
\endaligned
\end{equation}
Indeed,  we  take
$$
\aligned
\widetilde{\phi}(y)=\left(\|\phi\|_{**}+
\|h\|_{*}
\right)\psi(y),
\,
\qquad
\forall
\,\,y\in\overline{\widetilde{\Omega}}_t=
 \overline{B}_{2d/\varepsilon_0}(q')\setminus\left[\bigcup_{i=1}^mB_{R_1\gamma_i}(\xi'_i)
\cup B_{R_1\rho_0 v_0/\varepsilon_0}(q')\right].
\endaligned
$$
For
$y\in B_{2d/\varepsilon_0}(q')\setminus\big[\bigcup_{i=1}^mB_{R_1\gamma_i}(\xi'_i)\cup B_{R_1\rho_0 v_0/\varepsilon_0}(q')\big]$,
$$
\aligned
\mathcal{L}\big(\widetilde{\phi}\pm\phi\big)(y)
\geq C_1\|h\|_{*}\left\{\left(\frac{\varepsilon_0}{\rho_0v_0}\right)^2\frac{1}
{\big|\frac{\varepsilon_0 y-q}{\rho_0v_0}\big|^{4+2\hat{\alpha}}}
+
\sum_{i=1}^m\frac{1}{\gamma_i^2}\frac{1}{\big|\frac{y-\xi'_i}{\gamma_i}\big|^{4+2\hat{\alpha}}}
+\varepsilon_0^2
\right\}\pm
h(y)
\geq
|h(y)|\pm h(y)\geq0,
\endaligned
$$
and for $y\in\bigcup_{i=1}^m\partial B_{R_1\gamma_i}(\xi'_i)
\cup\partial B_{R_1\rho_0 v_0/\varepsilon_0}(q')\cup\partial B_{2d/\varepsilon_0}(q')$,
$$
\aligned
\big(\widetilde{\phi}\pm\phi\big)(y)
\geq\|\phi\|_{**}
\pm\phi(y)\geq|\phi(y)|
\pm\phi(y)
\geq0.
\endaligned
$$
By the above maximum
principle   we obtain  that $|\phi|\leq\tilde{\phi}$
in
$\widetilde{\Omega}_t=
\Omega_t\setminus\big[\bigcup_{i=1}^mB_{R_1\gamma_i}(\xi'_i)
\cup B_{R_1\rho_0 v_0/\varepsilon_0}(q')\cup B^c_{2d/\varepsilon_0}(q')\big]$,
which implies that estimate (\ref{3.12}) holds.

We prove the lemma by contradiction.  Assume that there are  sequences of
parameters $t_n\rightarrow+\infty$,
points $\xi^n=(\xi_1^n,\ldots,\xi_m^n)\in\mathcal{O}_{t_n}(q)$,
functions $h_n$, $W_n$  and associated solutions $\phi_n$ of
equation (\ref{3.10}) with orthogonality conditions (\ref{3.14})
such that
\begin{equation}\label{3.16}
\aligned
\|\phi_n\|_{L^{\infty}(\Omega_{t_n})}=1
\,\,\,\,\ \quad\,\,\,\,
\textrm{but}
\
\,\,\,\,\quad\,\,\,\,
\|h_n\|_{*}\rightarrow0,
\,\,\quad\,\,\textrm{as}\,\,\,\,n\rightarrow+\infty.
\endaligned
\end{equation}
Consider
$$
\aligned
\widehat{\phi}^n_{q^c}(x)
=\phi_n\big(x/\varepsilon_0^n\big),
\,\,\qquad\,\,
\,\,\qquad\,\,
\widehat{h}^n_{q^c}(x)=h_n\big(x/\varepsilon_0^n\big)
\,\,\qquad\,\,\textrm{for all}\,\,\,x\in\Omega\setminus B_{2d}(q),
\endaligned
$$
and
$$
\aligned
\widehat{\phi}^n_q(z)=\phi_n\big((\rho_0^nv_0^nz+q)/\varepsilon_0^n\big),
\,\,\quad\,\,
\qquad
\,\,\quad\,\,
\widehat{h}^n_q(z)=h_n\big((\rho_0^nv_0^nz+q)/\varepsilon_0^n\big),
\endaligned
$$
and for all $i=1,\ldots,m$,
$$
\aligned
\widehat{\phi}^n_i(z)=\phi_n\big(\gamma_{i}^nz+(\xi^n_i)'\big),
\,\,\quad\,\,
\,\,\,\,\qquad\,\,\,\,
\,\,\quad\,\,
\widehat{h}^n_i(z)=h_n\big(\gamma_{i}^nz+(\xi^n_i)'\big),
\endaligned
$$
where
$\mu^n=\big(\mu^n_0,\mu^n_1,\ldots,\mu_m^n\big)$,
$\varepsilon_0^n=\exp\left\{-\frac12t_n\right\}$,
$\varepsilon^n_i=\exp\left\{-\frac12t_n\phi_1(\xi_i^n)\right\}$,
$\rho^n_0=(\varepsilon_0^n)^{\frac{1}{1+\alpha}}=\exp\big\{
-\frac1{2(1+\alpha)}t_n\big\}$,
$v^n_0=(\mu_0^n)^{\frac{1}{1+\alpha}}$ and
$\gamma^n_i=\frac1{\varepsilon_0^n}\varepsilon_i^n\mu^n_i
=\mu^n_i\exp\left\{
-\frac12t_n\big[\phi_1(\xi^n_i)-1\big]\right\}$.
Notice first that
$$
\aligned
h_n(y)=\big(-\Delta_{a(\varepsilon^n_0 y)}\phi_n-W_n\phi_n\big)\big|_{y=x/\varepsilon_0^n}
=(\varepsilon^n_0)^{2}\left[
-\Delta_{a(x)}\widehat{\phi}^n_{q^c}-(\varepsilon^n_0)^{-2}\widehat{W}_n\widehat{\phi}^n_{q^c}
\right](x),
\endaligned
$$
where
$$
\aligned
\widehat{W}_n(x)=W_n(x/\varepsilon_0^n).
\endaligned
$$
Applying the expansion of $W_n$ in (\ref{2.32}), we find  that
 $\widehat{\phi}^n_{q^c}(x)$ satisfies
$$
\aligned
\left\{\aligned
&-\Delta_{a(x)}\widehat{\phi}^n_{q^c}(x)+
O\left(\frac{e^{-t_n\phi_1(x)}}
{\,|x-q|^{4+2\alpha}\,}\prod_{i=1}^m\frac{1}{\,|x-\xi_i^n|^4\,}
\right)
\widehat{\phi}^n_{q^c}(x)
=\left(\frac{1}{\varepsilon^n_0}\right)^2\widehat{h}^n_{q^c}(x)
\,\,\quad\textrm{in}\,\quad
\Omega\setminus B_{2d}(q),\\[1mm]
&\widehat{\phi}^n_{q^c}(x)=0
\qquad\qquad\qquad\qquad\qquad\qquad
\qquad\,\,\,\qquad\qquad\qquad\qquad\qquad
\qquad\quad\quad
\textrm{on}\,\qquad
\partial\Omega.
\endaligned
\right.
\endaligned
$$
Thanks to   the definition of the $\|\cdot\|_{*}$-norm in (\ref{3.7}),
we find that
$\big(\frac{1}{\varepsilon^n_0}\big)^2\big|\widehat{h}^n_{q^c}(x)\big|\leq C\|h_n\|_{*}\rightarrow 0$
uniformly in $\Omega\setminus B_{2d}(q)$. By elliptic estimates,
$\widehat{\phi}^n_{q^c}$
converges uniformly in $\Omega\setminus B_{2d}(q)$ to a trivial  solution $\widehat{\phi}^{\infty}_{q^c}$,
namely $\widehat{\phi}^{\infty}_{q^c}\equiv0$ in $\Omega\setminus B_{2d}(q)$.

Next, we observe that
$$
\aligned
h_n(y)=\big(-\Delta_{a(\varepsilon^n_0 y)}\phi_n-W_n\phi_n\big)\left|_{y=\large\frac{\rho_0^nv_0^nz+q}{\varepsilon_0^n}}
=\left(\frac{\varepsilon_0^n}{\rho_0^nv_0^n}\right)^2\left[
-\Delta_{\widehat{a}_n(z)}\widehat{\phi}_q^n
-\left(\frac{\rho_0^nv_0^n}{\varepsilon_0^n}\right)^2\widehat{W}_n\widehat{\phi}_q^n
\right](z),
\right.
\endaligned
$$
where
$$
\aligned
\widehat{a}_n(z)=a\big(\rho_0^nv_0^nz+q\big),
\qquad\qquad\quad
\widehat{W}_n(z)=W_n\big((\rho_0^nv_0^nz+q)/\varepsilon_0^n\big).
\endaligned
$$
Using the expansion of $W_n$ in (\ref{2.30}),
we find that $\widehat{\phi}^n_q(z)$ solves
$$
\aligned
-\Delta_{\widehat{a}_n(z)}\widehat{\phi}^n_q(z)
-\frac{8(1+\alpha)^2|z|^{2\alpha}}{(1+|z|^{2(1+\alpha)})^2}
\Big[1
+O\left(|\rho^n_0v^n_0z|^{\sigma}\right)
+
o\left(1\right)
\Big]
\widehat{\phi}^n_q(z)=\left(\frac{\rho^n_0v^n_0}{\varepsilon^n_0}\right)^2\widehat{h}^n_q(z)
\endaligned
$$
for any $z\in B_{R_0+2}(0)$.
Owing to  the definition of the $\|\cdot\|_{*}$-norm in (\ref{3.7}), we have that
 for any $\theta\in\big(1, -1/\hat{\alpha}\big)$,
$\big(\frac{\rho^n_0v^n_0}{\varepsilon^n_0}\big)^2\widehat{h}^n_q\rightarrow 0$ in
$L^{\theta}\big(B_{R_0+2}(0)\big)$.
Since
$\frac{8(1+\alpha)^2|z|^{2\alpha}}{(1+|z|^{2(1+\alpha)})^2}$
is bounded in $L^{\theta}\big(B_{R_0+2}(0)\big)$,
standard elliptic regularity  implies that
$\widehat{\phi}^n_q$
converges uniformly over
compact subsets near the origin to a
bounded solution $\widehat{\phi}^{\infty}_q$ of equation
\begin{equation}\label{3.3}
\aligned
\Delta
\phi+\frac{8(1+\alpha)^2|z|^{2\alpha}}{(1+|z|^{2(1+\alpha)})^2}\phi=0
\,\quad\textrm{in}\,\,\,\mathbb{R}^2,
\endaligned
\end{equation}
which satisfies
\begin{equation}\label{3.17}
\aligned
\int_{\mathbb{R}^2}\chi \mathcal{Z}_q\widehat{\phi}_q^{\infty}=0.
\endaligned
\end{equation}
By the result of \cite{E,CY,EPW}, $\widehat{\phi}^{\infty}_q$
is proportional to $\mathcal{Z}_q$.
Since $\int_{\mathbb{R}^2}\chi \mathcal{Z}_q^2>0$, by
(\ref{3.17}) we find that $\widehat{\phi}^{\infty}_q\equiv0$ in $B_{R_1}(0)$.

Finally,  for each $i\in\{1,\ldots,m\}$, we  get
$$
\aligned
h_n(y)=\big(-\Delta_{a(\varepsilon^n_0 y)}\phi_n-W_n\phi_n\big)\big|_{y=\gamma_{i}^n z+(\xi^n_i)'}
=(\gamma_i^n)^{-2}\left[
-\Delta_{\widehat{a}_n(z)}\widehat{\phi}_i^n-(\gamma_i^n)^{2}\widehat{W}_n\widehat{\phi}_i^n
\right](z),
\endaligned
$$
where
$$
\aligned
\widehat{a}_n(z)=a\big(\varepsilon_0^n\gamma_{i}^nz+\xi^n_i\big),
\qquad\qquad\qquad
\widehat{W}_n(z)=W_n\big(\gamma_{i}^nz+(\xi^n_i)'\big).
\endaligned
$$
Employing  the expansion of $W_n$ in (\ref{2.28}) and
 elliptic regularity,  we have  that for each $i\in\{1,\ldots,m\}$,
$\widehat{\phi}^n_i$
converges uniformly over
compact subsets near the origin to a bounded solution $\widehat{\phi}^{\infty}_i$ of equation
\begin{equation}\label{3.4}
\aligned
\Delta
\phi+\frac{8}{(1+|z|^2)^2}\phi=0
\,\quad\textrm{in}\,\,\,\mathbb{R}^2,
\endaligned
\end{equation}
which satisfies
\begin{equation}\label{3.18}
\aligned
\int_{\mathbb{R}^2}\chi \mathcal{Z}_j\widehat{\phi}_i^{\infty}=0
\quad\,\,\,\textrm{for}\,\,\,\,j=0,\,1,\,2.
\endaligned
\end{equation}
Thus by the result of  \cite{BP,CL}, $\widehat{\phi}^{\infty}_i$ must be
a linear combination of $\mathcal{Z}_j$, $j=0,1,2$.
But $\int_{\mathbb{R}^2}\chi \mathcal{Z}_j^2>0$
and $\int_{\mathbb{R}^2}\chi \mathcal{Z}_j\mathcal{Z}_l=0$ for $j\neq l$.
Hence (\ref{3.18}) implies  $\widehat{\phi}^{\infty}_i\equiv0$ in $B_{R_1}(0)$.
As a result, by definition (\ref{3.11}) we obtain
$\lim_{n\rightarrow+\infty}\|\phi_n\|_{**}=0$. But
(\ref{3.12}) and (\ref{3.16}) tell us
$\liminf_{n\rightarrow+\infty}\|\phi_n\|_{**}>0$,
which is a contradiction.
\end{proof}

\vspace{1mm}
\vspace{1mm}
\vspace{1mm}

{\bf Step 3:} Proving  an a priori estimate for solutions to
(\ref{3.10}) that satisfy orthogonality conditions with respect to
$Z_{ij}$ for $j=1,2$ only.

\vspace{1mm}
\vspace{1mm}
\vspace{1mm}

\noindent{\bf Lemma 3.3.}\,\,{\it For $t$  large enough, if
$\phi$ solves {\upshape (\ref{3.10})} and satisfies
\begin{equation}\label{3.19}
\aligned
\int_{\Omega_t}\chi_iZ_{ij}\phi=0\,\,\,\,\,
\,\,\,\,\forall\,\,i=1,\ldots,m,\,\,j=1,2,
\endaligned
\end{equation}
then
\begin{equation}\label{3.20}
\aligned
\|\phi\|_{L^{\infty}(\Omega_t)}\leq Ct\, \|h\|_{*},
\endaligned
\end{equation}
where $C>0$ is independent of $t$.}

\vspace{1mm}
\vspace{1mm}

\begin{proof}
Let $R>R_0+1$ be a large but fixed number. We consider the functions
\begin{equation}\label{3.21}
\aligned
\widehat{Z}_{q}(y)=Z_{q}(y)-\frac{\varepsilon_0}{\rho_0v_0}
+a_{q}G(\varepsilon_0 y,q),
\,\qquad\qquad\,
\widehat{Z}_{i0}(y)=Z_{i0}(y)-\frac1{\gamma_i}
+a_{i0}G(\varepsilon_0 y,\xi_i),
\endaligned
\end{equation}
where
\begin{equation}\label{3.22}
\aligned
a_{q}=\frac{\varepsilon_0}{\rho_0v_0\big[H(q,q)-4\log(\rho_0v_0 R)\big]},
\,\qquad\quad\qquad\,
a_{i0}=\frac1{\gamma_i\big[H(\xi_i,\xi_i)-4\log(\varepsilon_0\gamma_i R)\big]}.
\endaligned
\end{equation}
From estimates (\ref{2.14})-(\ref{2.15}), we obtain
\begin{equation}\label{3.23}
\aligned
C_1|\log\varepsilon_0|\leq-\log(\rho_0v_0 R)
\leq C_2|\log\varepsilon_0|,
\qquad\qquad
C_1|\log\varepsilon_i|\leq-\log(\varepsilon_0\gamma_i R)
\leq C_2|\log\varepsilon_i|,
\endaligned
\end{equation}
and
\begin{equation}\label{3.24}
\aligned
\widehat{Z}_{q}(y)=O\left(
\frac{\varepsilon_0G(\varepsilon_0 y,q)}{\rho_0v_0|\log\varepsilon_0|}
\right),
\,\,\,\qquad\quad\qquad\,
\widehat{Z}_{i0}(y)=O\left(
\frac{G(\varepsilon_0 y,\xi_i)}{\gamma_i|\log\varepsilon_i|}
\right).
\endaligned
\end{equation}
Let   $\eta_1$ and $\eta_2$ be  radial smooth cut-off functions in $\mathbb{R}^2$ such that
$$
\aligned
&0\leq\eta_1\leq1;
\qquad
\,\,\,\,\,\,\,\eta_1\equiv1\,\,\textrm{in}\,B_R(0);\,\,\,\,\,\,\,
\qquad\,\,\,
\eta_1\equiv0\,\,\textrm{in}\,\mathbb{R}^2\setminus B_{R+1}(0);\\[1mm]
&0\leq\eta_2\leq1;
\qquad
\,\,\,\,\,\,\,\eta_2\equiv1\,\,\textrm{in}\,B_{3d}(0);\,\,\,
\qquad\,\,\,\,\,\,\,\eta_2\equiv0\,\,\textrm{in}\,\mathbb{R}^2\setminus B_{6d}(0).
\endaligned
$$
Without loss of generality we assume that
$d>0$ is  a  small but fixed number independent of $t$ such that
$B_{9d}(q)\subset\Omega$.
Set
\begin{equation}\label{3.25}
\aligned
\eta_{q1}(y)=
\eta_1\left(\frac{\big|\varepsilon_0y-q\big|}{\rho_0 v_0}\right),
\,\,\qquad\qquad\,\,
\eta_{i1}(y)=
\eta_1\left(\frac{\big|y-\xi_i'\big|}{\gamma_i}\right),
\endaligned
\end{equation}
and
\begin{equation}\label{3.26}
\aligned
\eta_{q2}(y)=
\eta_2\left(\varepsilon_0\big|y-q'\big|\right),
\,\,\,\qquad\,\,\qquad\,\,\,
\eta_{i2}(y)=
\eta_2\left(\varepsilon_0\big|y-\xi'_i\big|\right),
\endaligned
\end{equation}
and  define the two test functions
\begin{equation}\label{3.27}
\aligned
\widetilde{Z}_{q}=\eta_{q1}Z_{q}+(1-\eta_{q1})\eta_{q2}\widehat{Z}_{q},
\,\,\,\qquad\,\,\qquad\,\,\,
\widetilde{Z}_{i0}=\eta_{i1}Z_{i0}+(1-\eta_{i1})\eta_{i2}\widehat{Z}_{i0}.
\endaligned
\end{equation}

Given  $\phi$ satisfying  (\ref{3.10}) and   (\ref{3.19}),
we modify it so that the extra orthogonality conditions with
respect to $Z_{q}$ and $Z_{i0}$'s hold. We set
\begin{equation}\label{3.28}
\aligned
\widetilde{\phi}=\phi+d_q\widetilde{Z}_{q}+\sum\limits_{i=1}^{m}d_i\widetilde{Z}_{i0}
+\sum_{i=1}^m\sum\limits_{j=1}^{2}e_{ij}\chi_iZ_{ij},
\endaligned
\end{equation}
and adjust the coefficients $d_q$,
$d_i$ and $e_{ij}$ such that $\widetilde{\phi}$
 satisfies the orthogonality condition
\begin{equation}\label{3.15}
\aligned
\int_{\Omega_t}\chi_qZ_q\widetilde{\phi}=0
\,\qquad\,
\textrm{and}
\,\qquad\,\int_{\Omega_t}\chi_iZ_{ij}\widetilde{\phi}=0
\,\,\ \ \,\,\,
\textrm{for all}\,\,\,i=1,\ldots,m,\,\,j=0,1,2.
\endaligned
\end{equation}
Then
\begin{equation}\label{3.35}
\aligned
\mathcal{L}(\widetilde{\phi})=h+d_q\mathcal{L}(\widetilde{Z}_{q})+\sum\limits_{i=1}^{m}d_i\mathcal{L}(\widetilde{Z}_{i0})
+\sum_{i=1}^m\sum_{j=1}^2e_{ij}\mathcal{L}(\chi_iZ_{ij})
\,\quad\,\textrm{in}\,\,\,\,\,\Omega_t.
\endaligned
\end{equation}
If   (\ref{3.15}) holds, the previous lemma  allows us to conclude
\begin{equation}\label{3.36}
\aligned
\|\widetilde{\phi}\|_{L^{\infty}(\Omega_t)}\leq
C\left[\|h\|_{*}+|d_q|\big\|\mathcal{L}(\widetilde{Z}_{q})\big\|_{*}
+\sum\limits_{i=1}^{m}|d_i|\big\|\mathcal{L}(\widetilde{Z}_{i0})\big\|_{*}
+\sum_{i=1}^m\sum_{j=1}^2|e_{ij}|\big\|\mathcal{L}(\chi_i Z_{ij})\big\|_{*}
\right].
\endaligned
\end{equation}
Furthermore, using the definition of $\widetilde{\phi}$ again and the fact that
\begin{equation}\label{3.37}
\aligned
\big\|\widetilde{Z}_{q}\big\|_{L^{\infty}(\Omega_t)}\leq\frac{C\varepsilon_0}{\rho_0v_0},
\quad\qquad\quad\quad
\big\|\widetilde{Z}_{i0}\big\|_{L^{\infty}(\Omega_t)}\leq\frac{C}{\gamma_i},
\quad\qquad\quad\quad
\big\|\chi_iZ_{ij}\big\|_{L^{\infty}(\Omega_t)}\leq\frac{C}{\gamma_i},
\endaligned
\end{equation}
estimate (\ref{3.20}) is  a  direct consequence of the following two claims:

\vspace{1mm}
\vspace{1mm}
\vspace{1mm}
\vspace{1mm}

\noindent{\bf Claim 1.}\,\,{\it
The coefficients $d_q$,
$d_i$ and $e_{ij}$ are well defined and
\begin{equation}\label{3.32}
\aligned
\big\|\mathcal{L}(\widetilde{Z}_{q})\big\|_{*}\leq
\frac{C\varepsilon_0\log t}{\rho_0v_0|\log\varepsilon_0|},
\endaligned
\end{equation}
and
\begin{equation}\label{3.33}
\aligned
\big\|\mathcal{L}(\chi_iZ_{ij})\big\|_{*}\leq\frac{C}{\gamma_i},
\,\quad\qquad\qquad\,\,\,\quad\,
\big\|\mathcal{L}(\widetilde{Z}_{i0})\big\|_{*}\leq
\frac{C\log t}{\gamma_i|\log\varepsilon_i|}.
\endaligned
\end{equation}
}

\noindent{\bf Claim 2.}\,\,{\it
The following bounds hold:
\begin{equation}\label{3.34}
\aligned
|d_q|\leq C\frac{\rho_0v_0|\log\varepsilon_0|}{\varepsilon_0}\|h\|_{*},
\,\quad\quad\,\,\quad\,\,\quad
|d_i|\leq C\gamma_i|\log\varepsilon_i|\|h\|_{*},
\,\quad\quad\,\,\quad\,\,\quad
|e_{ij}|\leq C\gamma_i\log t\,\|h\|_{*}.
\endaligned
\end{equation}
}

\vspace{1mm}
\vspace{1mm}
\vspace{1mm}

\noindent{\bf Proof of Claim 1.}
First, we find $d_q$,
$d_i$ and $e_{ij}$.
Testing  (\ref{3.28}) against $\chi_iZ_{ij}$ and using
 orthogonality condition (\ref{3.15})  for $j=1,2$ and the fact that
$\chi_i\chi_k\equiv0$ if $i\neq k$, we readily find
\begin{equation}\label{3.29}
\aligned
e_{ij}=\left(-d_q\int_{\Omega_t}\chi_iZ_{ij}\widetilde{Z}_{q}
-\sum_{k\neq i}^md_k\int_{\Omega_t}\chi_iZ_{ij}\widetilde{Z}_{k0}
\right)\left/\int_{\Omega_t}\chi^2_iZ^2_{ij},
\,\,\quad\,
i=1,\ldots,m,\,\,j=1,2.
\right.
\endaligned
\end{equation}
Note that $\int_{\Omega_t}\chi^2_iZ^2_{ij}=c>0$ for all $i$, $j$, and
$$
\aligned
\int_{\Omega_t}\chi_iZ_{ij}\widetilde{Z}_{q}=
O\left(\frac{\varepsilon_0\gamma_i\log t}{\rho_0v_0|\log\varepsilon_0|}
\right),
\,\,\,\qquad\,\,\,
\int_{\Omega_t}\chi_iZ_{ij}\widetilde{Z}_{k0}=
O\left(\frac{\gamma_i\log t}{\gamma_k|\log\varepsilon_k|}
\right),\,
\quad\,\,\,k\neq i.
\endaligned
$$
Then
\begin{equation}\label{3.30}
\aligned
|e_{ij}|\leq C\left(|d_q|\frac{\varepsilon_0\gamma_i\log t}{\rho_0v_0|\log\varepsilon_0|}
+\sum_{k\neq i}^m|d_k|\frac{\gamma_i\log t}{\gamma_k|\log\varepsilon_k|}\right).
\endaligned
\end{equation}
We only need  to consider $d_q$  and $d_i$. Testing  (\ref{3.28}) against $\chi_qZ_q$ and $\chi_kZ_{k0}$,
respectively, and using
the orthogonality conditions in (\ref{3.15})
for $q$ and $j=0$, we obtain a system of $\mathcal{D}=(d_q,d_1,\ldots,d_m)$,
\begin{equation}\label{3.31}
\aligned
d_q\int_{\Omega_t}\chi_qZ_q\widetilde{Z}_{q}
+\sum_{i=1}^md_i\int_{\Omega_t}\chi_qZ_q\widetilde{Z}_{i0}=&
-\int_{\Omega_t}\chi_qZ_q\phi,
\\[1mm]
d_q\int_{\Omega_t}\chi_kZ_{k0}\widetilde{Z}_{q}
+\sum_{i=1}^md_i\int_{\Omega_t}\chi_kZ_{k0}\widetilde{Z}_{i0}=&
-\int_{\Omega_t}\chi_kZ_{k0}\phi,
\,\quad\quad\forall\,\,\, k=1,\ldots,m.
\endaligned
\end{equation}
But
$$
\aligned
\int_{\Omega_t}\chi_qZ_q\widetilde{Z}_{q}=
\int_{\Omega_t}\chi_qZ_q^2=C_1>0,
\,\,\,\qquad\qquad\,\,\,
\int_{\Omega_t}\chi_qZ_q\widetilde{Z}_{i0}=
O\left(\frac{\rho_0v_0\log t}{\varepsilon_0\gamma_i|\log\varepsilon_i|}
\right),
\endaligned
$$
and
$$
\aligned
\int_{\Omega_t}\chi_kZ_{k0}\widetilde{Z}_{q}=
O\left(
\frac{\varepsilon_0\gamma_k\log t}{\rho_0v_0|\log\varepsilon_0|}
\right),
\qquad
\int_{\Omega_t}\chi_kZ_{k0}\widetilde{Z}_{k0}=C_2>0,
\qquad
\int_{\Omega_t}\chi_kZ_{k0}\widetilde{Z}_{i0}=
O\left(\frac{\gamma_k\log t}{\gamma_i|\log\varepsilon_i|}
\right),
\quad
i\neq k.
\endaligned
$$
Let us denote the system (\ref{3.31}) as $\mathcal{M}\mathcal{D}=\mathcal{S}$,
where $\mathcal{S}$ is the right-hand member
and $\mathcal{M}$ is the coefficient matrix.
By the
above estimates,
$\mathcal{P}^{-1}\mathcal{M}\mathcal{P}$ is diagonally dominant and then invertible, where
$\mathcal{P}=\diag\left(\rho_0 v_0\big/\varepsilon_0,\,\gamma_1,\,\ldots,\,\gamma_m\right)$. Hence $\mathcal{M}$ is also invertible and
$\mathcal{D}=(d_q,d_1,\ldots,d_m)$ is well defined.

Let us prove  now  inequality  (\ref{3.32}).
Consider four regions
$$
\aligned
\Omega_1=\left\{\,\left|\frac{\varepsilon_0y-q}{\rho_0 v_0}\right|\leq R\right\},
\quad\quad\qquad\quad&\qquad\quad
\Omega_2=\left\{R<\left|\frac{\varepsilon_0y-q}{\rho_0 v_0}\right|\leq R+1\right\},\\[1mm]
\Omega_3=\left\{R+1<\left|\frac{\varepsilon_0y-q}{\rho_0 v_0}\right|\leq\frac{3d}{\rho_0v_0}\right\},
\quad\qquad&\qquad\quad
\Omega_4=\left\{\frac{3d}{\rho_0v_0}<\left|\frac{\varepsilon_0y-q}{\rho_0 v_0}\right|\leq
\frac{6d}{\rho_0v_0}\right\}.
\quad\quad\,\,\,\,\,\,
\endaligned
$$
Notice first that (\ref{3.2}) and (\ref{3.5}) imply
\begin{eqnarray}\label{3.38}
\Delta_{a(\varepsilon_0 y)}Z_{q}+\left(\frac{\varepsilon_0}{\rho_0v_0}\right)^2
\frac{8(1+\alpha)^2\big|\frac{\varepsilon_0 y-q}{\rho_0 v_0}\big|^{2\alpha}}
{\big(1+\big|\frac{\varepsilon_0 y-q}{\rho_0 v_0}\big|^{2(1+\alpha)}\big)^2}Z_{q}
=\varepsilon_0\nabla\log a(\varepsilon_0 y)\nabla Z_{q}
\qquad\qquad\qquad\qquad\qquad\qquad\,\,\,\,\qquad
&&\nonumber\\
=O\left(\varepsilon_0\left(\frac{\varepsilon_0}{\rho_0v_0}\right)^2
\left|\frac{\varepsilon_0 y-q}{\rho_0 v_0}\right|^{2\alpha}
\left[1+\frac{|\varepsilon_0 y-q|}{\rho_0 v_0}\right]^{-3-4\alpha}\right),
&&
\end{eqnarray}
and then, in $\Omega_1\cup\Omega_2$, by (\ref{2.30}),
\begin{eqnarray}\label{3.38a}
\mathcal{L}(Z_{q})=\left[
-\Delta_{a(\varepsilon_0 y)}Z_{q}-\left(\frac{\varepsilon_0}{\rho_0v_0}\right)^2
\frac{8(1+\alpha)^2\big|\frac{\varepsilon_0 y-q}{\rho_0 v_0}\big|^{2\alpha}}
{\big(1+\big|\frac{\varepsilon_0 y-q}{\rho_0 v_0}\big|^{2(1+\alpha)}\big)^2}Z_{q}\right]
+\left[
\left(\frac{\varepsilon_0}{\rho_0v_0}\right)^2
\frac{8(1+\alpha)^2\big|\frac{\varepsilon_0 y-q}{\rho_0 v_0}\big|^{2\alpha}}
{\big(1+\big|\frac{\varepsilon_0 y-q}{\rho_0 v_0}\big|^{2(1+\alpha)}\big)^2}
-
W\right]Z_{q}
\qquad\,
&&\nonumber\\
=
\left(\frac{\varepsilon_0}{\rho_0v_0}\right)^3
\frac{8(1+\alpha)^2\big|\frac{\varepsilon_0 y-q}{\rho_0v_0}\big|^{2\alpha}}{\big(1+\big|\frac{\varepsilon_0 y-q}{\rho_0v_0}\big|^{2(1+\alpha)}\big)^2}
O\left(|\varepsilon_0y-q|^{\sigma}\right)
+
O\left(\varepsilon_0\left(\frac{\varepsilon_0}{\rho_0v_0}\right)^2
\left|\frac{\varepsilon_0 y-q}{\rho_0 v_0}\right|^{2\alpha}
\left[1+\frac{|\varepsilon_0 y-q|}{\rho_0 v_0}\right]^{-3-4\alpha}\right).
\,\,
&&
\end{eqnarray}
Hence in  $\Omega_1$,
\begin{equation}\label{3.39}
\aligned
\mathcal{L}(\widetilde{Z}_{q})=\mathcal{L}(Z_{q})=
\left(\frac{\varepsilon_0}{\rho_0v_0}\right)^3
\left|\frac{\varepsilon_0 y-q}{\rho_0v_0}\right|^{2\alpha}\Big[
O\left(\varepsilon_0^\sigma|y-q'|^{\sigma}\right)
+
O\left(\rho_0v_0\right)\Big].
\endaligned
\end{equation}
In $\Omega_2$,
\begin{eqnarray}\label{3.40}
\mathcal{L}(\widetilde{Z}_{q})=\eta_{q1}\mathcal{L}(Z_{q})+
(1-\eta_{q1})\mathcal{L}(\widehat{Z}_{q})
-2\nabla\eta_{q1}\nabla(Z_{q}-\widehat{Z}_{q})-(Z_{q}-\widehat{Z}_{q})\Delta_{a(\varepsilon_0 y)} \eta_{q1}
\,\,\quad\,\,&&\nonumber\\[1mm]
=\mathcal{L}(Z_{q})+(1-\eta_{q1})W(Z_{q}-\widehat{Z}_{q})
-2\nabla\eta_{q1}\nabla(Z_{q}-\widehat{Z}_{q})-(Z_{q}-\widehat{Z}_{q})\Delta_{a(\varepsilon_0 y)}\eta_{q1}.
&&
\end{eqnarray}
Note that
\begin{equation}\label{3.41}
\aligned
Z_{q}-\widehat{Z}_{q}=\frac{\varepsilon_0}{\rho_0 v_0}
-a_{q}G(\varepsilon_0 y,q)=
\frac{\varepsilon_0}{\rho_0v_0\big[H(q,q)-4\log(\rho_0v_0 R)\big]}
\left[
4\log\frac{|\varepsilon_0y-q|}{R\rho_0v_0}+O\left(
\varepsilon_0^\sigma|y-q'|^\sigma\right)
\right],
\endaligned
\end{equation}
and then, in $\Omega_2$,
\begin{equation}\label{3.42}
\aligned
|Z_{q}-\widehat{Z}_{q}|=O\left(
\frac{\varepsilon_0}{R\rho_0v_0|\log\varepsilon_0|}
\right),
\,\,\quad\quad\,\quad\,\,\quad\quad\,
|\nabla\big(Z_{q}-\widehat{Z}_{q}\big)|=O\left(
\frac{\varepsilon_0^2}{R\rho_0^2v_0^2|\log\varepsilon_0|}
\right).
\endaligned
\end{equation}
Moreover $|\nabla\eta_{q1}|=O\big(\varepsilon_0/(\rho_0v_0)\big)$ and
$|\Delta_{a(\varepsilon_0 y)} \eta_{q1}|=O\big(\varepsilon_0^2/(\rho_0^2v_0^2)\big)$.
These, together with (\ref{2.30}) and (\ref{3.38a}), gives
that in $\Omega_2$,
\begin{equation}\label{3.43}
\aligned
\mathcal{L}(\widetilde{Z}_{q})
=O\left(\frac{\varepsilon_0^3}{R\rho_0^3v_0^3|\log\varepsilon_0|}
\right).
\endaligned
\end{equation}
In $\Omega_3$, by (\ref{3.21}), (\ref{3.27}) and  (\ref{3.38}),
\begin{eqnarray*}
\mathcal{L}(\widetilde{Z}_{q})=\mathcal{L}(\widehat{Z}_{q})
=\mathcal{L}(Z_{q})-\mathcal{L}(Z_{q}-\widehat{Z}_{q})
\qquad\quad\quad\quad\quad\quad\quad\quad\quad\quad\quad\quad\,\,\,\,\,
&&\nonumber\\[1.6mm]
\equiv \mathcal{A}_1+\mathcal{A}_2+O\left(\varepsilon_0\left(\frac{\varepsilon_0}{\rho_0v_0}\right)^2
\left[1+\frac{|\varepsilon_0 y-q|}{\rho_0 v_0}\right]^{-3-2\alpha}\right),
&&
\end{eqnarray*}
where
\begin{eqnarray*}
\mathcal{A}_1=\left[\left(\frac{\varepsilon_0}{\rho_0v_0}\right)^2\frac{8(1+\alpha)^2\big|\frac{\varepsilon_0 y-q}{\rho_0v_0}\big|^{2\alpha}}{\,\big(1+\big|\frac{\varepsilon_0 y-q}{\rho_0v_0}\big|^{2(1+\alpha)}\big)^2\,}
-W\right]Z_{q}
\,\,\qquad\,\,
\textrm{and}
\,\,\qquad\,\,
\mathcal{A}_2=
W\left[\frac{\varepsilon_0}{\rho_0v_0}
-a_{q}G(\varepsilon_0 y,q)\right].
&&
\end{eqnarray*}
For the estimates of these two terms, we split $\Omega_3$  into some subregions:
$$
\aligned
&\Omega_{q}=\left\{\,R+1<\left|\frac{\varepsilon_0y-q}{\rho_0 v_0}\right|\leq
\frac{1}{\rho_0v_0 t^{2\beta}}\,\right\},\\
\Omega_{3,k}=\Omega_3\,\bigcap\,&\big\{\,|y-\xi'_k|\leq1/(\varepsilon_0 t^{2\beta})\,\big\}
\,\,\,\,\quad\ \,\,\textrm{and}\ \ \quad\,\,
\widetilde{\Omega}_3=\Omega_3\setminus
\left[\bigcup_{k=1}^m\Omega_{3,k}\cup \Omega_{q}\right].
\endaligned
$$
From (\ref{2.30}) and (\ref{2.32}), we obtain
$$
\aligned
\mathcal{A}_1=\left\{
\aligned
&
\left(\frac{\varepsilon_0}{\rho_0v_0}\right)^3\frac{8(1+\alpha)^2\big|\frac{\varepsilon_0 y-q}{\rho_0v_0}\big|^{2\alpha}}{\,\big(1+\big|\frac{\varepsilon_0 y-q}{\rho_0v_0}\big|^{2(1+\alpha)}\big)^2\,}
\left[\,O\left(\varepsilon_0^\sigma|y-q'|^{\sigma}\right)
+O\left(\rho_0^\sigma v_0^\sigma\right)+\sum_{j=1}^mO\left(\varepsilon_j^\sigma\mu_j^\sigma\right)\right]
\,\,\ \,\textrm{in}\,\,\,\,\Omega_{q},\\[2mm]
&
O\left(\frac{\varepsilon_0^5}{\rho_0v_0}\mu_0^2t^{4\beta(2+\alpha)}
\right)
+O\left(\frac{\varepsilon_0^3}{\rho_0v_0}e^{-t\phi_1(\varepsilon_0 y)}t^{4\beta(2m+2+\alpha)}
\right)
\quad\qquad\qquad\qquad\quad
\qquad\,\,
\textrm{in}
\,\,\,\,\,
\widetilde{\Omega}_{3}.
\endaligned
\right.
\endaligned
$$
Moreover, by (\ref{3.23}) and  (\ref{3.41}),
$$
\aligned
\mathcal{A}_2=\left\{
\aligned
&\left(\frac{\varepsilon_0}{\rho_0v_0}\right)^3\frac{8(1+\alpha)^2\big|\frac{\varepsilon_0 y-q}{\rho_0v_0}\big|^{2\alpha}}{\,\big(1+\big|\frac{\varepsilon_0 y-q}{\rho_0v_0}\big|^{2(1+\alpha)}\big)^2\,}
O\left(\frac{\log |\varepsilon_0y-q|
-\log(R\rho_0v_0)+
\varepsilon_0^\sigma|y-q'|^\sigma}{|\log\varepsilon_0|}\right)
\,\,\,\,\,\,\textrm{in}\,\,\,\Omega_{q},\\[2mm]
&O\left(\frac{\varepsilon_0^3}{\rho_0v_0}t^{4\beta(2m+2+\alpha)}e^{-t\phi_1(\varepsilon_0 y)}
\right)
\qquad\qquad\qquad\qquad\qquad
\qquad\qquad\qquad\qquad
\,\,\,\qquad
\textrm{in}\,\,\,\,\widetilde{\Omega}_{3}.
\endaligned
\right.
\endaligned
$$
Then in $\Omega_{q}\cup \widetilde{\Omega}_{3}$,
\begin{equation}\label{3.44}
\aligned
\mathcal{L}(\widetilde{Z}_{q})=
\mathcal{L}(\widehat{Z}_{q})=\left(\frac{\varepsilon_0}{\rho_0v_0}\right)^3\frac{8(1+\alpha)^2\big|\frac{\varepsilon_0 y-q}{\rho_0v_0}\big|^{2\alpha}}{\,\big(1+\big|\frac{\varepsilon_0 y-q}{\rho_0v_0}\big|^{2(1+\alpha)}\big)^2\,}
O\left(\frac{\log |\varepsilon_0y-q|
-\log(R\rho_0v_0)}{|\log\varepsilon_0|}\right).
\endaligned
\end{equation}
In $\Omega_{3,k}$ with $k=1,\ldots,m$, by (\ref{2.28}),  (\ref{3.24}) and  (\ref{3.38}),
\begin{eqnarray}\label{3.45}
\mathcal{L}(\widetilde{Z}_{q})=
\left(\frac{\varepsilon_0}{\rho_0v_0}\right)^2\frac{8(1+\alpha)^2\big|\frac{\varepsilon_0 y-q}{\rho_0v_0}\big|^{2\alpha}}{
\big(1+\big|\frac{\varepsilon_0 y-q}{\rho_0v_0}\big|^{2(1+\alpha)}\big)^2}Z_{q}
-\left[
\Delta_{a(\varepsilon_0 y)} Z_{q}+\left(\frac{\varepsilon_0}{\rho_0v_0}\right)^2\frac{8(1+\alpha)^2\big|\frac{\varepsilon_0 y-q}{\rho_0v_0}\big|^{2\alpha}}{\big(1+\big|\frac{\varepsilon_0 y-q}{\rho_0v_0}\big|^{2(1+\alpha)}\big)^2}Z_{q}
\right]-W\widehat{Z}_{q}
&&\nonumber
\\[1.5mm]
=\left(\frac{\varepsilon_0}{\rho_0v_0}\right)^2\frac{8(1+\alpha)^2\big|\frac{\varepsilon_0 y-q}{\rho_0v_0}\big|^{2\alpha}}{
\big(1+\big|\frac{\varepsilon_0 y-q}{\rho_0v_0}\big|^{2(1+\alpha)}\big)^2}Z_{q}
+O\left(
\frac{1}{\gamma_{k}^2}\frac{8}{\big(1+\big|\frac{y-\xi'_k}{\gamma_{k}}\big|^2\big)^2}
\cdot\frac{\varepsilon_0G(\varepsilon_0 y,q)}{\rho_0v_0|\log\varepsilon_0|}
\right)
\qquad\qquad\qquad\qquad\,
&&\nonumber
\\[1.5mm]
=O\left(\frac{1}{\gamma_{k}^2}\frac{8}{\big(1+\big|\frac{y-\xi'_k}{\gamma_{k}}\big|^2\big)^2}
\cdot\frac{\varepsilon_0\log t}{\rho_0v_0|\log\varepsilon_0|}
\right).
\qquad\qquad\qquad\qquad\qquad
\qquad\qquad\qquad\qquad\qquad
\qquad\quad\,\,\,&&
\end{eqnarray}
In $\Omega_4$, thanks to
$|\nabla\eta_{q2}|=O\left(\varepsilon_0\right)$,
$|\Delta_{a(\varepsilon_0 y)}  \eta_{q2}|=O\left(\varepsilon_0^2\right)$,
\begin{equation}\label{3.48}
\aligned
|\widehat{Z}_{q}|=O\left(
\frac{\varepsilon_0}{\rho_0v_0|\log\varepsilon_0|}
\right)\,\,\qquad\qquad\,
\textrm{and}
\,\,\qquad\qquad\,
|\nabla\widehat{Z}_{q}|=O\left(
\frac{\varepsilon_0^2}{\rho_0v_0|\log\varepsilon_0|}
\right),
\endaligned
\end{equation}
by (\ref{3.38})  we get
\begin{eqnarray}\label{3.46}
\mathcal{L}(\widetilde{Z}_{q})=-\eta_{q2}\Delta_{a(\varepsilon_0 y)}  Z_{q}-\eta_{q2}W\widehat{Z}_{q}
-2\nabla\eta_{q2}\nabla \widehat{Z}_{q}-\widehat{Z}_{q}\Delta_{a(\varepsilon_0 y)} \eta_{q2}
\qquad\qquad\qquad\qquad\qquad\qquad\qquad\qquad\qquad
&&\nonumber\\
=-\eta_{q2}W\widehat{Z}_{q}
+
O\left(\left(\frac{\varepsilon_0}{\rho_0v_0}\right)^3
\left(\frac{|\varepsilon_0 y-q|}{\rho_0 v_0}\right)^{-4-2\alpha}
+
\varepsilon_0\left(\frac{\varepsilon_0}{\rho_0v_0}\right)^2
\left(\frac{|\varepsilon_0 y-q|}{\rho_0 v_0}\right)^{-3-2\alpha}
+
\frac{\varepsilon_0^3}{\rho_0v_0|\log\varepsilon_0|}
\right).
&&
\end{eqnarray}
From the previous choice of the number
$d$ we get that  for any $y\in\Omega_4$ and any $k=1,\ldots,m$,
$$
\aligned
|y-\xi'_k|\geq|y-q'|-|q'-\xi_k'|\geq\frac{3d}{\varepsilon_0}
-\frac{d}{\varepsilon_0}=\frac{2d}{\varepsilon_0}>\frac{1}{\varepsilon_0 t^{2\beta}}.
\endaligned
$$
This combined with (\ref{2.32}) gives
\begin{equation}\label{3.47}
\aligned
W=O\left(\varepsilon_0^2e^{-t\phi_1(\varepsilon_0 y)}\right)
\,\,\quad\,\,\textrm{in}\,\,\,\,\Omega_4.
\endaligned
\end{equation}
Hence in this region,
\begin{equation}\label{3.49}
\aligned
\mathcal{L}(\widetilde{Z}_{q})
=O\left(
\frac{\varepsilon_0^3}{\rho_0v_0|\log\varepsilon_0|}
\right).
\endaligned
\end{equation}
Combing   (\ref{3.39}), (\ref{3.43}), (\ref{3.44}), (\ref{3.45}) and (\ref{3.49})
with
the definition of the $\|\cdot\|_{*}$-norm in (\ref{3.7}),
we conclude
$$
\aligned
\big\|\mathcal{L}(\widetilde{Z}_{q})\big\|_{*}=O\left(
\frac{\varepsilon_0\log t}{\rho_0v_0|\log\varepsilon_0|}\right).
\endaligned
$$

The  inequalities in (\ref{3.33}) are easy to get as they are very
similar to the consideration of inequality (\ref{3.32}). We leave the details for
readers.

\vspace{1mm}
\vspace{1mm}
\vspace{1mm}

\noindent{\bf Proof of Claim 2.}
Testing equation (\ref{3.35}) against $a(\varepsilon_0 y)\widetilde{Z}_{q}$ and using estimates (\ref{3.36})-(\ref{3.37}),
we can derive that
$$
\aligned
d_q
\int_{\Omega_t}&a(\varepsilon_0 y)\widetilde{Z}_{q}\mathcal{L}(\widetilde{Z}_{q})
+
\sum_{k=1}^md_k
\int_{\Omega_t}a(\varepsilon_0 y)\widetilde{Z}_{q}\mathcal{L}(\widetilde{Z}_{k0})
\\
=&-\int_{\Omega_t}a(\varepsilon_0 y) h\widetilde{Z}_{q}
+\int_{\Omega_t}a(\varepsilon_0 y)\widetilde{\phi}\mathcal{L}(\widetilde{Z}_{q})-\sum_{k=1}^m\sum_{l=1}^2e_{kl}
\int_{\Omega_t}a(\varepsilon_0 y)\chi_kZ_{kl}\mathcal{L}(\widetilde{Z}_{q})
\\[1mm]
\leq&\frac{C\varepsilon_0}{\,\rho_0v_0\,}\|h\|_{*}
+C\big\|\mathcal{L}(\widetilde{Z}_{q})\big\|_{*}
\left(\|\widetilde{\phi}\|_{L^{\infty}(\Omega_t)}
+\sum_{k=1}^m\sum_{l=1}^2\frac{1}{\gamma_k}|e_{kl}|\right)
\\[1mm]
\leq&\frac{C\varepsilon_0}{\,\rho_0v_0\,}\|h\|_{*}
+
C\big\|\mathcal{L}(\widetilde{Z}_{q})\big\|_{*}
\left[\|h\|_{*}
+|d_q|\big\|\mathcal{L}(\widetilde{Z}_{q})\big\|_{*}
+
\sum\limits_{k=1}^{m}|d_k|\big\|\mathcal{L}(\widetilde{Z}_{k0})\big\|_{*}
+\sum_{k=1}^m\sum_{l=1}^2|e_{kl}|\left(\frac{1}{\gamma_k}
+\big\|\mathcal{L}(\chi_kZ_{kl})\big\|_{*}
\right)
\right],
\endaligned
$$
where we have  applied  that
\begin{equation}\label{3.61}
\aligned
\left(\frac{\varepsilon_0}{\rho_0v_0}\right)^2\int_{\Omega_t}
\frac{\big|\frac{\varepsilon_0 y-q}{\rho_0v_0}\big|^{2\alpha}}
{\,\big(1+\big|\frac{\varepsilon_0 y-q}{\rho_0v_0}\big|\big)^{4+2\hat{\alpha}+2\alpha}\,}
dy\leq C
\,\,\,\quad\,\,\textrm{and}\,\,\,\quad\,\,
\frac{1}{\gamma_i^2}
\int_{\Omega_t}\frac{1}{\big(1+\big|\frac{y-\xi'_i}{\gamma_i}\big|\big)^{4+2\hat{\alpha}}}
dy\leq C,
\,\,\,\,\,i=1,\ldots,m.
\endaligned
\end{equation}
But estimate (\ref{3.30}) and Claim $1$ imply
\begin{equation}\label{3.50}
\aligned
|d_q|\left|
\int_{\Omega_t}a(\varepsilon_0 y)\widetilde{Z}_{q}\mathcal{L}(\widetilde{Z}_{q})
\right|
\leq
\frac{C\varepsilon_0}{\rho_0v_0}\|h\|_{*}
+\frac{C\varepsilon_0\log^2 t}{\rho_0v_0|\log\varepsilon_0|}\left[
\frac{\varepsilon_0|d_q|}{\rho_0v_0|\log\varepsilon_0|}
+
\sum_{k=1}^m\frac{|d_k|}{\gamma_k|\log\varepsilon_k|}
\right]
+\sum_{k=1}^m\left|d_k
\int_{\Omega_t}a(\varepsilon_0 y)\widetilde{Z}_{k0}\mathcal{L}(\widetilde{Z}_{q})
\right|.
\endaligned
\end{equation}
Similarly, testing  (\ref{3.35}) against $a(\varepsilon_0 y)\widetilde{Z}_{i0}$ and using
(\ref{3.36}), (\ref{3.37}), (\ref{3.30}) and Claim $1$, we obtain
\begin{eqnarray}\label{3.51}
|d_i|\left|
\int_{\Omega_t}a(\varepsilon_0 y)\widetilde{Z}_{i0}\mathcal{L}(\widetilde{Z}_{i0})
\right|
\leq
\frac{C\|h\|_{*}}{\gamma_i}
+\frac{C\log^2t}{\gamma_i|\log\varepsilon_i|}
\left[
\frac{\varepsilon_0|d_q|}{\rho_0v_0|\log\varepsilon_0|}
+\sum_{k=1}^m\frac{|d_k|}{\gamma_k|\log\varepsilon_k|}
\right]
+
\left|d_q
\int_{\Omega_t}a(\varepsilon_0 y)\widetilde{Z}_{i0}\mathcal{L}(\widetilde{Z}_{q})
\right|
&&\nonumber\\[1mm]
+
\sum_{k\neq i}^m\left|d_k
\int_{\Omega_t}a(\varepsilon_0 y)\widetilde{Z}_{k0}\mathcal{L}(\widetilde{Z}_{i0})
\right|.
\qquad\qquad\qquad\qquad\qquad\,
\qquad\qquad\qquad\qquad\quad\,\,\,\,
&&
\end{eqnarray}
To achieve the estimates  of $d_q$, $d_i$ and $e_{ij}$ in (\ref{3.34}),  we have the following claim.

\vspace{1mm}
\vspace{1mm}
\vspace{1mm}

\noindent{\bf Claim 3.}\,\,{\it
If  $d$ is sufficiently small,
but $R$  is  sufficiently large, then we have that
for any $i,k=1,\ldots,m$ with $i\neq k$,
\begin{equation}\label{3.53}
\aligned
\int_{\Omega_t}a(\varepsilon_0 y)\widetilde{Z}_{i0}\mathcal{L}(\widetilde{Z}_{i0})
=\frac{2\pi a(\xi_i)}{\gamma_i^2|\log\varepsilon_i|}\left[1+O\left(\frac1{R^2}\right)\right],
\qquad
\qquad
\int_{\Omega_t}a(\varepsilon_0 y)\widetilde{Z}_{k0}\mathcal{L}(\widetilde{Z}_{i0})
=O\left(\frac{\log^2t}{\gamma_i\gamma_k|\log\varepsilon_i||\log\varepsilon_k|}
\right),
\endaligned
\end{equation}
and
\begin{equation}\label{3.54}
\aligned
\int_{\Omega_t}a(\varepsilon_0 y)\widetilde{Z}_{q}\mathcal{L}(\widetilde{Z}_{q})
=\frac{2\pi(1+\alpha)a(q)}{|\log\varepsilon_0|}\left(\frac{\varepsilon_0}{\rho_0v_0}\right)^2
\left[1+O\left(\frac1{R^{2(1+\alpha)}}\right)\right],
\quad
\int_{\Omega_t}a(\varepsilon_0 y)\widetilde{Z}_{k0}\mathcal{L}(\widetilde{Z}_{q})
=O\left(\frac{\varepsilon_0\log^2 t}{\rho_0v_0\gamma_k|\log\varepsilon_0||\log\varepsilon_k|}\right).
\endaligned
\end{equation}
}

\vspace{1mm}

Indeed, once  Claim $3$ is valid, then  replacing (\ref{3.53}) and (\ref{3.54})
in (\ref{3.51}) and (\ref{3.50}), respectively,  we conclude
\begin{equation}\label{3.55}
\aligned
\frac{\varepsilon_0|d_q|}{\rho_0v_0|\log\varepsilon_0|}
\leq
C\|h\|_{*}
+\frac{C\log^2 t}{\,|\log\varepsilon_0|\,}\left(
\frac{\varepsilon_0|d_q|}{\rho_0v_0|\log\varepsilon_0|}
+
\sum_{k=1}^m\frac{|d_k|}{\gamma_k|\log\varepsilon_k|}
\right),
\endaligned
\end{equation}
and  for any $i=1,\ldots,m$,
\begin{equation}\label{3.56}
\aligned
\frac{|d_i|}{\gamma_i|\log\varepsilon_i|}
\leq
C\|h\|_{*}
+\frac{C\log^2t}{|\log\varepsilon_i|}
\left(
\frac{\varepsilon_0|d_q|}{\rho_0v_0|\log\varepsilon_0|}
+\sum_{k=1}^m\frac{|d_k|}{\gamma_k|\log\varepsilon_k|}
\right).
\endaligned
\end{equation}
As a result,  using linear algebra arguments for (\ref{3.55})-(\ref{3.56}), by (\ref{2.7})  we can prove  Claim 2 for
$d_q$ and $d_i$,  and then
complete the proof by inequality (\ref{3.30}).

\vspace{1mm}
\vspace{1mm}
\vspace{1mm}

\noindent{\bf Proof of Claim 3.}
Let us first establish the validity of the two expansions in (\ref{3.54}).
We decompose
\begin{equation*}\label{3.52}
\aligned
\int_{\Omega_t}a(\varepsilon_0 y)\widetilde{Z}_{q}\mathcal{L}(\widetilde{Z}_{q})
=\sum_{l=1}^4\int_{\Omega_l}a(\varepsilon_0 y)\widetilde{Z}_{q}\mathcal{L}(\widetilde{Z}_{q})
\equiv\sum_{l=1}^4I_{l}.
\endaligned
\end{equation*}
From (\ref{3.5}) and (\ref{3.39}) we obtain
\begin{eqnarray*}
I_1=\int_{\Omega_1}a(\varepsilon_0 y)Z_{q}\mathcal{L}(Z_{q})
=\int_{\Omega_1}
\left(\frac{\varepsilon_0}{\rho_0v_0}\right)^4
\left|\frac{\varepsilon_0 y-q}{\rho_0v_0}\right|^{2\alpha}\Big[
O\left(\varepsilon_0^\sigma|y-q'|^{\sigma}\right)
+
O\left(\rho_0v_0\right)\Big]
=\left(\frac{\varepsilon_0}{\rho_0v_0}\right)^2
O\left(\rho_0^\sigma v_0^\sigma\right)
\end{eqnarray*}
From (\ref{3.24}), (\ref{3.44}) and (\ref{3.45}), we deduce
\begin{eqnarray*}
I_3=\int_{\Omega_{q}\cup \widetilde{\Omega}_{3}}a(\varepsilon_0 y)\widehat{Z}_{q}\mathcal{L}(\widehat{Z}_{q})
+\sum_{k=1}^m\int_{\Omega_{3,k}}a(\varepsilon_0 y)\widehat{Z}_{q}\mathcal{L}(\widehat{Z}_{q})
\qquad\qquad\qquad\qquad\quad\,\,\,
\qquad\qquad\qquad\qquad\qquad
&&\nonumber\\[1mm]
=\left(\frac{\varepsilon_0}{\rho_0v_0}\right)^2
O
\left(\int_{R+1}^{3d/(\rho_0v_0)}\frac{r^{1+2\alpha}\log(r/R)}{(1+r^{2(1+\alpha)})^2}
\frac{\log(\rho_0v_0 r)}{|\log\varepsilon_0|^2}dr
+\sum_{k=1}^m\int_{0}^{1/(\varepsilon_0\gamma_kt^{2\beta})}
\frac{r}{(1+r^2)^2}
\frac{\log^2 t}{|\log\varepsilon_0|^2} dr
\right)
&&\nonumber\\[1.5mm]
=\left(\frac{\varepsilon_0}{\rho_0v_0}\right)^2
\left[O\left(\frac{1}{R^{2(1+\alpha)}|\log\varepsilon_0|}\right)
+O\left(\frac{\log^2t}{|\log\varepsilon_0|^2}\right)\right].
\qquad\qquad\qquad\qquad\,\,\,\,
\qquad\qquad\qquad\qquad\quad\quad
&&
\end{eqnarray*}
From (\ref{3.48}) and (\ref{3.49}), we derive that
\begin{equation*}
\aligned
I_4=\int_{\Omega_4}a(\varepsilon_0 y)\eta_{q2}\widehat{Z}_{q}\mathcal{L}(\widetilde{Z}_{q})
=\int_{\left\{\frac{3d}{\rho_0v_0}<\left|\frac{\varepsilon_0y-q}{\rho_0 v_0}\right|\leq
\frac{6d}{\rho_0v_0}\right\}}
O\left(
\frac{\varepsilon_0^4}{\rho_0^2v_0^2|\log\varepsilon_0|^2}
\right)
dy
=O\left(
\frac{\varepsilon_0^2}{\rho_0^2v_0^2|\log\varepsilon_0|^2}
\right).
\endaligned
\end{equation*}
As for $I_2$, by (\ref{3.40}) we  get
$$
\aligned
I_2=&-\int_{\Omega_2}a(\varepsilon_0 y)\widetilde{Z}_{q}(Z_{q}-\widehat{Z}_{q})\Delta_{a(\varepsilon_0 y)}\eta_{q1}
-2\int_{\Omega_2}a(\varepsilon_0 y)\widetilde{Z}_{q}\nabla\eta_{q1}\nabla(Z_{q}-\widehat{Z}_{q})
\\
&+\int_{\Omega_2}a(\varepsilon_0 y)\widetilde{Z}_{q}\big[
\mathcal{L}(Z_{q})+(1-\eta_{q1})W(Z_{q}-\widehat{Z}_{q})
\big].
\endaligned
$$
Integrating by parts the first term and using estimates (\ref{2.30}), (\ref{3.38}) and (\ref{3.42}) for the last term, we obtain
\begin{eqnarray*}
I_2=
-\int_{\Omega_2}a(\varepsilon_0 y)\widehat{Z}_{q}\nabla\eta_{q1}\nabla(Z_{q}-\widehat{Z}_{q})
+\int_{\Omega_2}a(\varepsilon_0 y)(Z_{q}-\widehat{Z}_{q})^2|\nabla\eta_{q1}|^2
\,\,\,
&&\nonumber\\[1.5mm]
+\int_{\Omega_2}a(\varepsilon_0 y)(Z_{q}-\widehat{Z}_{q})\nabla\eta_{q1}\nabla\widehat{Z}_{q}
+\left(\frac{\varepsilon_0}{\rho_0v_0}\right)^2
O\left(\frac{1}{R^{3+2\alpha}|\log\varepsilon_0|}\right)
&&\nonumber\\[1.5mm]
=I_{21}+I_{22}+I_{23}+\left(\frac{\varepsilon_0}{\rho_0v_0}\right)^2
O\left(\frac{1}{R^{3+2\alpha}|\log\varepsilon_0|}\right).
\qquad\qquad\qquad\quad\,
&&
\end{eqnarray*}
From  (\ref{3.2}), (\ref{3.5}), (\ref{3.25})  and (\ref{3.42}) we find
$|\nabla\eta_{q1}|=O\big(\frac{\varepsilon_0}{\rho_0v_0}\big)$
and $|\nabla\widehat{Z}_{q}|=O\big(\frac{\varepsilon_0^2}{R^{3+2\alpha}\rho_0^2v_0^2}\big)$ in $\Omega_2$.
Furthermore,
\begin{equation*}
\aligned
I_{22}=O\left(\frac{\varepsilon_0^2}{R\rho_0^2v_0^2|\log\varepsilon_0|^2}\right)
\quad\qquad\quad
\textrm{and}
\quad\qquad\quad
I_{23}=O\left(\frac{\varepsilon_0^2}{R^{3+2\alpha}\rho_0^2v_0^2|\log\varepsilon_0|}\right).
\endaligned
\end{equation*}
Since  $a(\varepsilon_0 y)=a(q)\big[1+O(\varepsilon_0 |y-q'|)\big]$
and $\widehat{Z}_{q}=Z_{q}\big[1+O\big(\frac{1}{R|\log\varepsilon_0|}\big)\big]$
in $\Omega_2$,
by   (\ref{2.7a}), (\ref{2.14}), (\ref{3.2}), (\ref{3.5}),  (\ref{3.25}) and (\ref{3.41})  we conclude
\begin{eqnarray*}
&&I_{21}=-\left(\frac{\varepsilon_0}{\rho_0v_0}\right)^3
\frac{1}{H(q,q)-4\log(\rho_0v_0 R)}
\int_{\left\{R<\left|\frac{\varepsilon_0y-q}{\rho_0 v_0}\right|\leq R+1\right\}}
\frac{a(\varepsilon_0 y)}{|y-q'|}
\mathcal{Z}_q\left(\frac{\varepsilon_0y-q}{\rho_0 v_0}\right)
\eta_1'\left(\frac{\big|\varepsilon_0y-q\big|}{\rho_0 v_0}\right)
\big(4+o(1)\big)dy\nonumber\\[1.5mm]
&&\,\quad\,\,=-\left(\frac{\varepsilon_0}{\rho_0v_0}\right)^2
\frac{8\pi a(q)}{H(q,q)-4\log(\rho_0v_0 R)}
\int_{R}^{R+1}
\eta_1'(r)\left[
1+O\left(\frac1{r^{2(1+\alpha)}}\right)
\right]dr\nonumber\\[2mm]
&&\,\quad\,\,=\frac{2\pi(1+\alpha)a(q)}{|\log\varepsilon_0|}\left(\frac{\varepsilon_0}{\rho_0v_0}\right)^2
\left[1+O\left(\frac1{R^{2(1+\alpha)}}\right)\right].
\end{eqnarray*}
Combining all these estimates, we have
that for
$R$ and $t$ large enough, and  $d$ small enough,
\begin{equation}\label{3.59}
\aligned
\int_{\Omega_t}a(\varepsilon_0 y)\widetilde{Z}_{q}\mathcal{L}(\widetilde{Z}_{q})
=\frac{2\pi(1+\alpha)a(q)}{|\log\varepsilon_0|}\left(\frac{\varepsilon_0}{\rho_0v_0}\right)^2
\left[\,1+O\left(\frac1{R^{2(1+\alpha)}}\right)\right].
\endaligned
\end{equation}

According to (\ref{3.50}), we only need to consider  $\int_{\Omega_t}a(\varepsilon_0 y)\widetilde{Z}_{k0}\mathcal{L}(\widetilde{Z}_{q})$
for all $k$.
By the above estimates of
$\mathcal{L}(\widetilde{Z}_{q})$  and $\widetilde{Z}_{k0}$, we have clearly
$$
\aligned
\int_{\Omega_1}a(\varepsilon_0 y)
\widetilde{Z}_{k0}\mathcal{L}(\widetilde{Z}_{q})
=O\left(\frac{\varepsilon_0
\rho_0^\sigma v_0^\sigma
\log t
}{\rho_0v_0\gamma_k|\log\varepsilon_k|}
\right),
\qquad\qquad\qquad
\int_{\Omega_2}a(\varepsilon_0 y)
\widetilde{Z}_{k0}\mathcal{L}(\widetilde{Z}_{q})
=O\left(\frac{\varepsilon_0\log t}
{\rho_0v_0\gamma_k|\log\varepsilon_0||\log\varepsilon_k|}\right),
\endaligned
$$
$$
\aligned
\int_{\Omega_4}a(\varepsilon_0 y)\widetilde{Z}_{k0}\mathcal{L}(\widetilde{Z}_{q})=
O\left(\frac{\varepsilon_0 }{\rho_0v_0\gamma_k|\log\varepsilon_0||\log\varepsilon_k|}\right),
\,\,\quad\qquad\,\,
\int_{\Omega_{q}\cup \widetilde{\Omega}_{3}}a(\varepsilon_0 y)\widetilde{Z}_{k0}\mathcal{L}(\widetilde{Z}_{q})
=O\left(\frac{\varepsilon_0\log t}{\rho_0v_0\gamma_k|\log\varepsilon_0||\log\varepsilon_k|}\right),
\endaligned
$$
and
$$
\aligned
\int_{\Omega_{3,l}}a(\varepsilon_0 y)\widetilde{Z}_{k0}\mathcal{L}(\widetilde{Z}_{q})
=O\left(\frac{\varepsilon_0\log^2 t}{\rho_0v_0\gamma_k|\log\varepsilon_0||\log\varepsilon_k|}\right)
\quad\quad\textrm{for all}\,\,\,l\neq k.
\endaligned
$$
It remains to  calculate  the integral over $\Omega_{3,k}$. Using (\ref{3.27}) and  an integration by parts,  we obtain
$$
\aligned
\int_{\Omega_{3,k}}a(\varepsilon_0 y)\widetilde{Z}_{k0}\mathcal{L}(\widetilde{Z}_{q})
=\int_{\Omega_{3,k}}a(\varepsilon_0 y)\widehat{Z}_{q}\mathcal{L}(\widetilde{Z}_{k0})
-\int_{\partial\Omega_{3,k}}a(\varepsilon_0 y)\widehat{Z}_{k0}\frac{\partial\widehat{Z}_{q}}{\partial\nu}
+\int_{\partial\Omega_{3,k}}a(\varepsilon_0 y)\widehat{Z}_{q}\frac{\partial\widehat{Z}_{k0}}{\partial\nu}.
\endaligned
$$
Let us decompose
$$
\aligned
\int_{\Omega_{3,k}}a(\varepsilon_0 y)\widehat{Z}_{q}\mathcal{L}(\widetilde{Z}_{k0})
=
\left(\int_{\big\{|y-\xi'_k|\leq\gamma_kR\big\}}
+
\int_{\big\{\gamma_kR<|y-\xi'_k|\leq\gamma_k(R+1)\big\}}
+
\int_{\big\{\gamma_k(R+1)<|y-\xi'_k|\leq
1/(\varepsilon_0t^{2\beta})\big\}}
\right)
a(\varepsilon_0 y)
\widehat{Z}_{q}\mathcal{L}(\widetilde{Z}_{k0}).
\endaligned
$$
In a straightforward but tedious way,  we can compute that for $|y-\xi'_k|\leq\gamma_kR$,
$$
\aligned
\mathcal{L}(\widetilde{Z}_{k0})=\mathcal{L}(Z_{k0})=
O\left(\rho_0^\sigma v_0^\sigma/\gamma_k^3\right)+\sum_{j=1}^mO\left(\varepsilon_j^\sigma\mu_j^\sigma/\gamma_k^3\right),
\endaligned
$$
for $\gamma_kR<|y-\xi'_k|\leq\gamma_k(R+1)$,
$$
\aligned
\mathcal{L}(\widetilde{Z}_{k0})
=O\left(\frac{1}{R\gamma_k^3|\log\varepsilon_k|}
\right),
\endaligned
$$
and for $\gamma_k(R+1)<|y-\xi'_k|\leq
1/(\varepsilon_0t^{2\beta})$,
$$
\aligned
\mathcal{L}(\widetilde{Z}_{k0})=
\mathcal{L}(\widehat{Z}_{k0})=O
\left(\frac{\log|y-\xi'_k|-
\log(R\gamma_k)}{\big(1+\big|\frac{y-\xi'_k}{\gamma_{k}}\big|^2\big)^2}\cdot
\frac{1}{\gamma_k^3|\log\varepsilon_k|}
\right).
\endaligned
$$
These combined the estimate of $\widehat{Z}_{q}$ in (\ref{3.24}) give
$$
\aligned
\int_{\Omega_{3,k}}
a(\varepsilon_0 y)\widehat{Z}_{q}\mathcal{L}(\widetilde{Z}_{k0})
=O\left(\frac{\varepsilon_0\log t}{\rho_0v_0\gamma_k|\log\varepsilon_0||\log\varepsilon_k|}\right).
\endaligned
$$
As on $\partial\Omega_{3,k}$, by (\ref{2.2}) and (\ref{3.24}) we know
$$
\aligned
\widehat{Z}_{q}
=O\left(\frac{\varepsilon_0\log t}{\rho_0v_0|\log\varepsilon_0|}\right),
\,\qquad\,\,\,\qquad\,
|\nabla\widehat{Z}_{q}|
=O\left(\frac{\varepsilon_0^2t^{\beta}}{\rho_0v_0|\log\varepsilon_0|}\right),
\endaligned
$$
and
$$
\aligned
\widehat{Z}_{k0}
=O\left(\frac{\log t}{\gamma_k|\log\varepsilon_k|}\right),
\,\qquad\,\,\,\qquad\,
|\nabla\widehat{Z}_{k0}|
=O\left(\frac{\varepsilon_0  t^{2\beta} }{\gamma_k|\log\varepsilon_k|}\right).
\endaligned
$$
Thus
$$
\aligned
\int_{\Omega_{3,k}}a(\varepsilon_0 y)\widetilde{Z}_{k0}\mathcal{L}(\widetilde{Z}_{q})
=O\left(\frac{\varepsilon_0\log t}{\rho_0v_0\gamma_k|\log\varepsilon_0||\log\varepsilon_k|}\right).
\endaligned
$$
By the above estimates, we readily  have
\begin{equation}\label{3.60}
\aligned
\int_{\Omega_t}a(\varepsilon_0 y)\widetilde{Z}_{k0}\mathcal{L}(\widetilde{Z}_{q})
=O\left(\frac{\varepsilon_0\log^2 t}{\rho_0v_0\gamma_k|\log\varepsilon_0||\log\varepsilon_k|}\right),
\,\quad\,k=1,\ldots,m.
\endaligned
\end{equation}

The  two expansions in (\ref{3.53}) are easy to get as they are very
similar to the above consideration for the two expansions in (\ref{3.54}). We leave the details for
readers.
\end{proof}

\vspace{1mm}
\vspace{1mm}
\vspace{1mm}

{\bf Step 4:}
Proof of Proposition 3.1. Let us first establish the validity of the a priori
estimate (\ref{3.8}).
From the previous lemma and the fact that $\|\chi_iZ_{ij}\|_{*}=O(\gamma_i)$, we give
\begin{equation}\label{3.65}
\aligned
\|\phi\|_{L^{\infty}(\Omega_t)}\leq Ct\left(
\|h\|_{*}+\sum\limits_{i=1}^m\sum\limits_{j=1}^2\gamma_i|c_{ij}|
\right),
\endaligned
\end{equation}
hence it is sufficient to estimate the size of the constants $c_{ij}$.
Let us consider the cut-off function $\eta_{i2}$ introduced in
 (\ref{3.26}).  Testing    (\ref{3.1}) against  $a(\varepsilon_0 y)\eta_{i2}Z_{ij}$,
$i=1,\ldots,m$ and $j=1,2$, we obtain
\begin{equation}\label{3.66}
\aligned
\int_{\Omega_t}
a(\varepsilon_0 y)
\phi \mathcal{L}(\eta_{i2}Z_{ij})=
\int_{\Omega_t}
a(\varepsilon_0 y)
h\eta_{i2}Z_{ij}+
\sum_{k=1}^m\sum_{l=1}^2c_{kl}
\int_{\Omega_t}\chi_kZ_{kl}\eta_{i2}Z_{ij}.
\endaligned
\end{equation}
For any  $i=1,\ldots,m$ and $j=1,2$,
$$
\aligned
\Delta_{a(\varepsilon_0 y)}Z_{ij}+\frac{1}{\gamma_{i}^2}\frac{8}{\big(1+\big|\frac{y-\xi'_i}{\gamma_{i}}\big|^2\big)^2}Z_{ij}
=\varepsilon_0\nabla\log a(\varepsilon_0 y)\nabla Z_{ij}
=O\left(\frac{\varepsilon_0}{\gamma_{i}^2}
\left[1+\frac{|y-\xi'_i|}{\gamma_{i}}\right]^{-2}\right).
\endaligned
$$
Then
$$
\aligned
\mathcal{L}(\eta_{i2}Z_{ij})=&\,
\eta_{i2}\mathcal{L}(Z_{ij})
-Z_{ij}\Delta_{a(\varepsilon_0 y)} \eta_{i2}
-2\nabla\eta_{i2}\nabla Z_{ij}\\
=&\left[\frac{1}{\gamma_{i}^2}\frac{8}{\big(1+\big|\frac{y-\xi'_i}{\gamma_{i}}\big|^2\big)^2}-W\right]\eta_{i2}Z_{ij}
-\eta_{i2}\left[
\Delta_{a(\varepsilon_0 y)}Z_{ij}+\frac{1}{\gamma_{i}^2}\frac{8}{\big(1+\big|\frac{y-\xi'_i}{\gamma_{i}}\big|^2\big)^2}Z_{ij}
\right]
+O\left(\varepsilon_0^3\right)\\
\equiv&\,\mathcal{B}_{ij}+O\left(\frac{\varepsilon_0}{\gamma_{i}^2}
\left[1+\frac{|y-\xi'_i|}{\gamma_{i}}\right]^{-2}\right)
+
O\left(\varepsilon_0^3\right),
\endaligned
$$
where
$$
\aligned
\mathcal{B}_{ij}
=\left[\frac{1}{\gamma_{i}^2}\frac{8}{\big(1+\big|\frac{y-\xi'_i}{\gamma_{i}}\big|^2\big)^2}-W\right]\eta_{i2}Z_{ij}.
\endaligned
$$
For the estimate of   $\mathcal{B}_{ij}$, we decompose  $\supp(\eta_{i2})$ to some subregions:
$$
\aligned
\widehat{\Omega}_{q}=\supp(\eta_{i2})\bigcap\big\{|y-q'|\leq
1/(\varepsilon_0 t^{2\beta})\big\},
\,\quad\quad\quad\,
\widehat{\Omega}_{k1}=\supp(\eta_{i2})\bigcap\big\{|y-\xi'_k|\leq
1/(\varepsilon_0 t^{2\beta})\big\},
\,\,\,\,k=1,\ldots,m,
\endaligned
$$
$$
\aligned
\widehat{\Omega}_{2}=\supp(\eta_{i2})\setminus\left[\,\bigcup_{k=1}^m\widehat{\Omega}_{k1}
\cup\widehat{\Omega}_{q}\right],
\endaligned
$$
where $\supp(\eta_{i2})=\left\{|y-\xi'_i|\leq
6d/\varepsilon_0\right\}$.
Notice that, by (\ref{2.2}),
\begin{equation}\label{3.67}
\aligned
|y-\xi'_i|\geq |\xi_i'-q'|-|y-q'|\geq
|\xi_i'-q'|-\frac{1}{\varepsilon_0 t^{2\beta}}
\geq
\frac{1}{\varepsilon_0 t^\beta}\left(1-\frac{1}{t^\beta}
\right)
>\frac{1}{2\varepsilon_0 t^\beta}
\endaligned
\end{equation}
uniformly in $\widehat{\Omega}_{q}$, and
\begin{equation}\label{3.68}
\aligned
|y-\xi'_i|\geq |\xi_i'-\xi_k'|-|y-\xi_k'|\geq
|\xi_i'-\xi_k'|-\frac{1}{\varepsilon_0 t^{2\beta}}
\geq
\frac{1}{\varepsilon_0 t^\beta}\left(1-\frac{1}{t^\beta}
\right)>\frac{1}{2\varepsilon_0 t^\beta}
\endaligned
\end{equation}
uniformly in $\widehat{\Omega}_{k1}$ with $k\neq i$.
From expansions  (\ref{2.30})-(\ref{2.32}) of $W$
and definition (\ref{3.6})  of $Z_{ij}$  we get,
in $\widehat{\Omega}_{i1}$,
$$
\aligned
\mathcal{B}_{ij}=
\frac{1}{\gamma_{i}^3}\frac{8}{\big(1+\big|\frac{y-\xi'_i}{\gamma_{i}}\big|^2\big)^{5/2}}
\left[O\left(\varepsilon_0^\sigma|y-\xi'_i|^{\sigma}\right)
+O\left(\rho_0^\sigma v_0^\sigma\right)+\sum_{j=1}^mO\left(\varepsilon_j^\sigma\mu_j^\sigma\right)
\right],
\endaligned
$$
and in $\widehat{\Omega}_{q}$,  by (\ref{3.67}),
$$
\aligned
\mathcal{B}_{ij}=
\left[
O\left(\frac{\gamma_i^2}{|y-\xi'_i|^4}\right)+
\left(\frac{\varepsilon_0}{\rho_0v_0}\right)^2
O\left(\frac{8(1+\alpha)^2\big|\frac{\varepsilon_0 y-q}{\rho_0v_0}\big|^{2\alpha}}{\,\big(1+\big|\frac{\varepsilon_0 y-q}{\rho_0v_0}\big|^{2(1+\alpha)}\big)^2\,}\right)
\right]O\left(\varepsilon_0 t^\beta
\right),
\endaligned
$$
and in $\widehat{\Omega}_{k1}$, $k\neq i$, by (\ref{3.68}),
$$
\aligned
\mathcal{B}_{ij}=\left[
O\left(\frac{\gamma_i^2}{|y-\xi'_i|^4}\right)+
O\left(\frac{1}{\gamma_{k}^2}\frac{8}{\big(1+\big|\frac{y-\xi'_k}{\gamma_{k}}\big|^2\big)^2}\right)
\right]O\left(\varepsilon_0 t^\beta
\right),
\endaligned
$$
and in $\widehat{\Omega}_{2}$,
$$
\aligned
\mathcal{B}_{ij}=
O\left(\varepsilon_0^3\varepsilon_i^2\mu_i^2t^{10\beta}\right)
+O\left(\varepsilon_0^3 t^{2\beta(4m+5+2\alpha)} e^{-t\phi_1(\varepsilon_0 y)}
\right).
\endaligned
$$
So,
$$
\aligned
\left|\int_{\Omega_t}
a(\varepsilon_0 y)\phi \mathcal{L}(\eta_{i2}Z_{ij})\right|
\leq C\frac1{\gamma_i}\max\big\{(\rho_0 v_0)^\sigma,(\varepsilon_1\mu_1)^\sigma,
\ldots,(\varepsilon_m\mu_m)^\sigma
\big\}\|\phi\|_{L^{\infty}(\Omega_t)}.
\endaligned
$$
On the other hand, since
$\|\eta_{i2}Z_{ij}\|_{L^{\infty}(\Omega_t)}\leq C\gamma_i^{-1}$, by (\ref{3.7}) and (\ref{3.61}) we know that
$$
\aligned
\int_{\Omega_t}a(\varepsilon_0 y)h\eta_{i2}Z_{ij}=O\left(\frac{\|h\|_{*}}{\gamma_i}\right).
\endaligned
$$
Now,  if $k= i$, by (\ref{3.2}), (\ref{3.5}) and (\ref{3.6}),
$$
\aligned
\int_{\Omega_t}\chi_kZ_{kl}\eta_{k2}Z_{kj}=
\int_{\mathbb{R}^2}\chi(|z|)\mathcal{Z}_{l}(z) \mathcal{Z}_{j}(z)dz=C\delta_{lj},
\endaligned
$$
while if $k\neq i$, by (\ref{3.68}),
$$
\aligned
\int_{\Omega_t}\chi_kZ_{kl}\eta_{i2}Z_{ij}=O\left(\gamma_k\varepsilon_0 t^\beta\right).
\endaligned
$$
Using the above estimates  in (\ref{3.66}),  we find
$$
\aligned
|c_{ij}|\leq C\left(\frac1{\gamma_i}\max\big\{(\rho_0 v_0)^\sigma,(\varepsilon_1\mu_1)^\sigma,\ldots,(\varepsilon_m\mu_m)^\sigma
\big\}\|\phi\|_{L^{\infty}(\Omega_t)}+
\frac{1}{\gamma_i}\|h\|_{*}+
\sum\limits_{k\neq i}^m\sum\limits_{l=1}^2\gamma_k\varepsilon_0 t^\beta|c_{kl}|
\right),
\endaligned
$$
and then
$$
\aligned
|c_{ij}|\leq C\frac1{\gamma_i}\Big(\max\big\{(\rho_0 v_0)^\sigma,(\varepsilon_1\mu_1)^\sigma,\ldots,(\varepsilon_m\mu_m)^\sigma
\big\}\|\phi\|_{L^{\infty}(\Omega_t)}+
\|h\|_{*}
\Big).
\endaligned
$$
Putting this estimate in (\ref{3.65}), we conclude the validity of  (\ref{3.8}).

Let us consider the Hilbert space
$$
\aligned
H_{\xi}=\left\{\phi\in H_0^1(\Omega_t)\left|\,
\int_{\Omega_t}\chi_iZ_{ij}\phi=0\,\,\,\,\,\,\forall\,\,\,i=1,\ldots,m,\,\,j=1,2\right.\right\}
\endaligned
$$
with the norm
$\|\phi\|_{H_\xi}=\|\nabla\phi\|_{L^{2}(\Omega_t)}$.
Equation (\ref{3.1}) is equivalent to find
$\phi\in H_\xi$, such that
$$
\aligned
\int_{\Omega_t}a(\varepsilon_0 y)\nabla\phi\nabla\psi-\int_{\Omega_t}a(\varepsilon_0 y)W\phi\psi
=\int_{\Omega_t}a(\varepsilon_0 y)h\psi,\,\,\quad\,\,\forall
\,\,\psi\in H_\xi.
\endaligned
$$
By Fredholm's alternative this is equivalent to the uniqueness of solutions to this
problem, which is  guaranteed by   (\ref{3.8}).
\end{proof}


\vspace{1mm}
\vspace{1mm}
\vspace{1mm}

\noindent{\bf Lemma 3.4.}\,\,{\it
For any integer $m\geq1$, the
operator $\mathcal{T}$ is differentiable  with
respect to the variables $\xi=(\xi_1,\ldots,\xi_m)$ in $\mathcal{O}_t(q)$,
precisely for any $k=1,\ldots,m$ and $l=1,2$,
\begin{equation}\label{3.74}
\aligned
\|\partial_{\xi'_{kl}}\mathcal{T}(h)\|_{L^{\infty}(\Omega_t)}\leq Ct^2\|h\|_{*}.
\endaligned
\end{equation}
}

To prove this lemma, we give the following estimate.

\vspace{1mm}
\vspace{1mm}
\vspace{1mm}
\vspace{1mm}

\noindent{\bf Lemma 3.5.}\,\,{\it
For any $0<\sigma<1$ and $\xi=(\xi_1,\ldots,\xi_m)\in\mathcal{O}_t(q)$, then we have
that for any $k=1,\ldots,m$ and $l=1,2$,
$$
\aligned
\partial_{\xi'_{kl}}H_i(x)=O\left((\varepsilon_{k}\mu_{k})^{\sigma}\right),
\endaligned
$$
uniformly in $\overline{\Omega}$,
where  $H_i$, $i=0,1,\ldots,m$, is defined as the solution of  equation {\upshape(\ref{2.9})}.
}

\vspace{1mm}
\vspace{1mm}

\begin{proof}
Differentiating   equation (\ref{2.9}) of $H_0$'s  with respect to  $\xi_{kl}$,  we  obtain
$$
\left\{\aligned
&-\Delta_{a}
\big(\partial_{\xi_{kl}}H_0\big)=\nabla\log a(x)\nabla\big(\partial_{\xi_{kl}}u_0\big)
\,\,\,\,
\textrm {in}\,\,\,\,\,\Omega,\\[2mm]
&\partial_{\xi_{kl}}H_0=-\partial_{\xi_{kl}}u_0
\,\qquad\qquad\qquad\qquad\quad\,
\textrm{on}\,\,\,\po,
\endaligned\right.
$$
where
$$
\aligned
\partial_{\xi_{kl}}u_0=2\partial_{\xi_{kl}}\log
\mu_0-
\frac{4\varepsilon_{0}^2\mu_{0}^2\partial_{\xi_{kl}}\log\mu_{0}}
{\varepsilon_{0}^2\mu_{0}^2+|x-q|^{2(1+\alpha)}}
\qquad
\textrm{and}
\qquad
\big|\nabla\big(\partial_{\xi_{kl}}u_0\big)\big|\leq
4\varepsilon_{0}^2\mu_{0}^2\big|\partial_{\xi_{kl}}\log\mu_{0}
\big|\frac{2(1+\alpha)|x-q|^{1+2\alpha}}
{(\varepsilon_{0}^2\mu_{0}^2+|x-q|^{2(1+\alpha)})^2}.
\endaligned
$$
Clearly,
$$
\aligned
\big\|\partial_{\xi_{kl}}u_0\big\|_{C^2(\partial\Omega)}
\leq C\big|\partial_{\xi_{kl}}\log\mu_0\big|\leq Ct^{\beta}.
\endaligned
$$
For any $\max\{1,\,2/(3+2\alpha)\}<p<2$,  by  the change of  variables $\rho_{0}v_{0}z=x-q$  we give
$$
\aligned
\int_{\Omega}\left|
\frac{
|x-q|^{1+2\alpha}}
{(\varepsilon_{0}^2\mu_{0}^2+|x-q|^{2(1+\alpha)})^2}
\right|^pdx=
\int_{\Omega_{\rho_{0}v_{0}}}
\frac{|z|^{(1+2\alpha)p}}
{(1+|z|^{2(1+\alpha)})^{2p}}(\rho_{0}v_{0})^{2-(3+2\alpha)p}
dz=O\left(\frac1{(\rho_{0}v_{0})^{(3+2\alpha)p-2}}
\right),
\\[0.1mm]
\endaligned
$$
where $\Omega_{\rho_{0}v_{0}}=\frac{1}{\rho_{0}v_{0}}(\Omega-\{q\})$.
Then
$$
\aligned
\left\|\nabla\log a(x)\nabla\big(\partial_{\xi_{kl}}u_0\big)\right\|_{L^p(\Omega)}
\leq C\left\|\nabla\big(\partial_{\xi_{kl}}u_0\big)\right\|_{L^p(\Omega)}
\leq C\frac{\varepsilon_{0}^2\mu_{0}^2\big|\partial_{\xi_{kl}}\log\mu_{0}
\big|}{(\rho_{0}v_{0})^{[(3+2\alpha)p-2]/p}}
\leq Ct^{\beta}(\rho_{0}v_{0})^{(2-p)/p}.
\endaligned
$$
Applying $L^p$ theory,
$$
\aligned
\big\|\partial_{\xi_{kl}}H_0\big\|_{W^{2,p}(\Omega)}\leq
C\left(\left\|\Delta_{a}
\big(\partial_{\xi_{kl}}H_0\big)\right\|_{L^p(\Omega)}
+\big\|\partial_{\xi_{kl}}H_0\big\|_{C^2(\partial\Omega)}\right)
\leq Ct^{\beta}.
\endaligned
$$
By Sobolev embedding, we find that for any $0<\gamma<2-(2/p)$,
$$
\aligned
\big\|\partial_{\xi_{kl}}H_0\big\|_{C^{\gamma}(\overline{\Omega})}
\leq C
t^{\beta}.
\endaligned
$$
Similarly, using the above arguments for each $H_i$, $i=1,\ldots,m$ again, we  readily
get
$$
\aligned
\big\|\partial_{\xi_{kl}}H_k\big\|_{C^{\gamma}(\overline{\Omega})}
\leq C (\varepsilon_{k}\mu_{k})^{(2-2p)/p}
\qquad\quad
\textrm{but}
\qquad\quad
\big\|\partial_{\xi_{kl}}H_i\big\|_{C^{\gamma}(\overline{\Omega})}
\leq C t^{\beta}
\qquad
\forall\,\,i\neq k.
\endaligned
$$
By (\ref{2.1}),   (\ref{2.7}),   (\ref{2.7a}),   (\ref{2.14}),  (\ref{2.15})
and  $\xi'=\xi/\varepsilon_0$, the lemma is then proven.
\end{proof}

\vspace{1mm}
\vspace{1mm}

\noindent{\bf Proof of Lemma 3.4.}
Differentiating (\ref{3.1}) with respect to $\xi'_{kl}$, formally
$Z=\partial_{\xi'_{kl}}\phi$ should satisfy
$$
\aligned
\left\{
\aligned
&\mathcal{L}(Z)=\phi\,\partial_{\xi'_{kl}}W+
\frac{1}{a(\varepsilon_0 y)}\sum\limits_{i=1}^m\sum\limits_{j=1}^2\left[c_{ij}\partial_{\xi'_{kl}}(\chi_iZ_{ij})+\widetilde{c}_{ij}\chi_iZ_{ij}\right]\,\,\quad
\ \textrm{in}\,\,\,\,\,\Omega_t,\\
&Z=0
\quad\qquad\qquad\qquad\quad
\,\quad\quad\quad\,\,\,\,
\qquad\qquad\qquad\quad
\qquad\qquad\qquad\,\,
\textrm{on}
\,\,\,\po_t,\\[1mm]
&\int_{\Omega_t}\chi_iZ_{ij}Z=-\int_{\Omega_t}\phi\partial_{\xi'_{kl}}(\chi_iZ_{ij})
\qquad\qquad\quad\quad\,
\forall\,\,\,i=1,\ldots,m,
\,\,j=1,2,
\endaligned
\right.
\endaligned
$$
where (still formally) $\widetilde{c}_{ij}=\partial_{\xi'_{kl}}(c_{ij})$.
We consider the constants $b_{ij}$ defined as
$$
\aligned
b_{ij}\int_{\Omega_t}\chi^2_i|Z_{ij}|^2=\int_{\Omega_t}\phi\,\partial_{\xi'_{kl}}(\chi_iZ_{ij}),
\endaligned
$$
and set
$$
\aligned
\widetilde{Z}=Z+\sum_{i=1}^m\sum_{j=1}^2b_{ij}\chi_iZ_{ij}.
\endaligned
$$
Then
$$
\aligned
\left\{
\aligned
&\mathcal{L}(\widetilde{Z})=f+\frac{1}{a(\varepsilon_0 y)}\sum\limits_{i=1}^m\sum\limits_{j=1}^2\widetilde{c}_{ij}\chi_iZ_{ij}
\quad\quad\,\,\,
\textrm{in}\,\,\,\,\,\,\,\Omega_t,\\
&\widetilde{Z}=0
\,\quad\qquad\quad\qquad\qquad\qquad\qquad\qquad
\,\,\,\,\,\textrm{on}
\,\,\,\,\,\po_t,\\[1mm]
&\int_{\Omega_t}\chi_iZ_{ij}\widetilde{Z}=0\,\,\,\,
\,\qquad\quad\,\forall\,\,i=1,\ldots,m,\,\,j=1,2,
\endaligned
\right.
\endaligned
$$
where
$$
\aligned
f=\phi\,\partial_{\xi'_{kl}}W
+\sum\limits_{i=1}^m\sum\limits_{j=1}^2
b_{ij}\mathcal{L}(\chi_iZ_{ij})
+\frac{1}{a(\varepsilon_0 y)}
\sum\limits_{i=1}^m
\sum\limits_{j=1}^2c_{ij}\partial_{\xi'_{kl}}(\chi_iZ_{ij}).
\endaligned
$$
The result of Proposition 3.1 implies  that this equation  has a unique solution
$\widetilde{Z}$ and $\widetilde{c}_{ij}$, and thus
$\partial_{\xi'_{kl}}\mathcal{T}(h)
=\mathcal{T}(f)-\sum_{i=1}^m\sum_{j=1}^2b_{ij}\chi_iZ_{ij}$
is well defined. Moreover,
\begin{equation}\label{3.75}
\aligned
\|\partial_{\xi'_{kl}}T(h)\|_{L^{\infty}(\Omega_t)}\leq\|\mathcal{T}(f)\|_{L^{\infty}(\Omega_t)}
+C\sum_{i=1}^m\sum_{j=1}^2\frac{1}{\gamma_i}|b_{ij}|
\leq C t\|f\|_{*}+C\sum_{i=1}^m\sum_{j=1}^2\frac{1}{\gamma_i}|b_{ij}|.
\endaligned
\end{equation}
To prove estimate (\ref{3.74}), we
 first estimate $\partial_{\xi'_{kl}}W$.
Notice that $\partial_{\xi'_{kl}}W=W\partial_{\xi'_{kl}}V$. Obviously,
by (\ref{2.30}),  (\ref{2.28}), (\ref{2.32}) and (\ref{3.7}) we  have that
$\|W\|_{*}=O\left(1\right)$. Moreover, thanks to Lemma $3.5$,
  by (\ref{2.4}),
(\ref{2.8}), (\ref{2.14}), (\ref{2.15}),  (\ref{3.2}) and (\ref{3.6})
we can directly check that
\begin{equation}\label{3.91}
\aligned
\partial_{\xi'_{kl}}V(y)=Z_{kl}(y)+O\left((\varepsilon_{k}\mu_{k})^{\sigma}\right).
\endaligned
\end{equation}
This, together with the fact that $\frac{1}{\gamma_k}\leq C$ uniformly on $t$, immediately implies
$$
\aligned
\|\partial_{\xi'_{kl}}V\|_{L^{\infty}(\Omega_t)}=O\left(1\right)
\,\,\qquad\,\,\textrm{and}\,\,\qquad\,\,\|\partial_{\xi'_{kl}}W\|_{*}=O\left(1\right).
\endaligned
$$
Next, by definitions (\ref{3.2}) and (\ref{3.6}),
a straightforward but tedious computation shows that
$$
\aligned
\|\partial_{\xi'_{kl}}(\chi_iZ_{ij})\|_{*}=\left\{
\aligned
&O\left(|\partial_{\xi'_{kl}}\gamma_i|\right)
\,\,\,\,\,\textrm{if}\,\,\,\,i\neq k,\\
&O\left(1\right)\,\,\quad
\,\,\,\,\,\,\quad\,\,\textrm{if}
\,\,\,\,i=k,
\endaligned
\right.
\,\,\qquad\qquad\,\,
|b_{ij}|=\left\{
\aligned
&O\left(|\partial_{\xi'_{kl}}\gamma_i|\right)\|\phi\|_{L^{\infty}(\Omega_t)}
\,\,\,\,\,\textrm{if}\,\,\,\,i\neq k,\\
&O\left(1\right)\|\phi\|_{L^{\infty}(\Omega_t)}\,\,\quad
\,\,\,\,\,\quad\,\,\,\textrm{if}
\,\,\,\,i=k.
\endaligned
\right.
\endaligned
$$
Furthermore, using  (\ref{3.8}),
(\ref{3.33})
and  the fact that
$\big|\partial_{\xi'_{kl}}\gamma_i\big|=O(\varepsilon_0^\sigma)$ uniformly on $t$,   we find
$$
\aligned
\|f\|_{*}\leq Ct\|h\|_{*}
\,\,\,\,\qquad\,\,\,\,
\textrm{and}
\,\,\,\,\qquad\,\,\,\,
|b_{ij}|\leq Ct\|h\|_{*}.
\endaligned
$$
Substituting these estimates into (\ref{3.75}), we then prove  (\ref{3.74}).
\qquad\qquad\qquad\qquad\qquad\qquad\qquad
\qquad\qquad\qquad\qquad\qquad\qquad$\square$

\vspace{1mm}
\vspace{1mm}
\vspace{1mm}

\section{The  intermediate  nonlinear problem}
In order to solve problem (\ref{2.34}), we shall  solve first the
intermediate nonlinear  problem: for any integer $m\geq1$
and
any points
$\xi=(\xi_1,\ldots,\xi_m)\in\mathcal{O}_t(q)$,
we find a function $\phi$ and scalars $c_{ij}$,
$i=1,\ldots,m$, $j=1,2$,  such that
\begin{equation}\label{4.1}
\aligned
\left\{
\aligned
&\mathcal{L}(\phi)=-\Delta_{a(\varepsilon_0 y)}\phi-W\phi=E+N(\phi)+\frac1{a(\varepsilon_0 y)}\sum\limits_{i=1}^m\sum\limits_{j=1}^2c_{ij}\chi_iZ_{ij}
\quad
\textrm{in}\,\,\,\,\,\Omega_t,\\
&\,\phi=0
\quad\quad\quad\quad\qquad\qquad\qquad\quad\ \,
\qquad\qquad\qquad\qquad\,
\quad\quad\quad\quad\qquad
\textrm{on}
\,\,\,\po_t,\\[1mm]
&\int_{\Omega_t}\chi_iZ_{ij}\phi=0
\,\,\qquad\qquad\qquad\qquad\,\,\,
\qquad\qquad\quad\ \,
\forall\,\,\,i=1,\ldots,m,
\,\,j=1,2,
\endaligned
\right.
\endaligned
\end{equation}
where  $W$ is as in (\ref{2.30}),  (\ref{2.28})  and (\ref{2.32}),
and $E$, $N(\phi)$ are given by
(\ref{2.21}) and (\ref{2.35}), respectively.

\vspace{1mm}
\vspace{1mm}
\vspace{1mm}
\vspace{1mm}

\noindent{\bf Proposition 4.1.}\,\,{\it
Let $m$ be a positive   integer and
 $0<\sigma<\min\{1/2,\,1-1/(2\beta),\,2(1+\alpha)\}$.
Then there exist constants $t_m>1$ and $C>0$ such
that for any $t>t_m$ and any points
$\xi=(\xi_1,\ldots,\xi_m)\in\mathcal{O}_t(q)$,
problem {\upshape(\ref{4.1})} admits
a unique solution
$\phi\in L^{\infty}(\Omega_t)$, and
scalars $c_{ij}\in\mathbb{R}$,
$i=1,\ldots,m$, $j=1,2$, such that
\begin{eqnarray}\label{4.2}
\quad\quad
\|\phi\|_{L^{\infty}(\Omega_t)}\leq C t
\max\big\{(\rho_0v_0)^{\min\{\sigma,2(\alpha-\hat{\alpha})\}},\,(\varepsilon_1\mu_1)^{\sigma},\,\ldots,\,(\varepsilon_m\mu_m)^{\sigma},\,
\|e^{-\frac{1}{2}t\phi_1}\|_{L^{\infty}(\Omega_t)}\big\}.
\end{eqnarray}
Furthermore, the map
$\xi'\mapsto\phi_{\xi'}\in C(\overline{\Omega}_t)$
is $C^1$, precisely for any $k=1,\ldots,m$ and  $l=1,2$,
\begin{eqnarray}\label{4.3}
\quad\quad
\|\partial_{\xi'_{kl}}\phi\|_{L^{\infty}(\Omega_t)}
\leq C t^{2}\max\big\{(\rho_0v_0)^{\min\{\sigma,2(\alpha-\hat{\alpha})\}},\,(\varepsilon_1\mu_1)^{\sigma},\,\ldots,\,(\varepsilon_m\mu_m)^{\sigma},\,
\|e^{-\frac{1}{2}t\phi_1}\|_{L^{\infty}(\Omega_t)}\big\},
\end{eqnarray}
where $\xi':=(\xi'_1,\ldots,\xi'_m)=(\frac1{\varepsilon_0}\xi_1,\ldots,\frac1{\varepsilon_0}\xi_m)$.
}

\vspace{1mm}
\vspace{1mm}

\begin{proof}
Proposition 3.1 and Lemma 3.4 allow us to apply
the contraction mapping  theorem and the implicit function theorem
to find a unique solution for problem (\ref{4.1})
satisfying (\ref{4.2})-(\ref{4.3}).
Since it is a standard procedure, we shall not present the details here, see Lemmas $4.1$-$4.2$ in \cite{DKM}
for a similar proof.
But we just mention that
$\|N(\phi)\|_{*}\leq C\|\phi\|^2_{L^{\infty}(\Omega_t)}$,
$\|E\|_{*}\leq C\max\{(\rho_0v_0)^{\min\{\sigma,2(\alpha-\hat{\alpha})\}},(\varepsilon_1\mu_1)^{\sigma},
\ldots,(\varepsilon_m\mu_m)^{\sigma},
\|e^{-\frac{1}{2}t\phi_1}\|_{L^{\infty}(\Omega_t)}\}$ and $\|\partial_{\xi'_{kl}}E\|_{*}
\leq C\max\{(\rho_0v_0)^{\min\{\sigma,2(\alpha-\hat{\alpha})\}},(\varepsilon_1\mu_1)^{\sigma},\ldots,
(\varepsilon_m\mu_m)^{\sigma},
\|e^{-\frac{1}{2}t\phi_1}\|_{L^{\infty}(\Omega_t)}\}$ due to (\ref{2.31})-(\ref{2.33})  and (\ref{3.91}).
\end{proof}

\vspace{1mm}

\section{The reduced problem: A maximization procedure }
In this section we study a maximization problem. Let
us consider the energy functional  associated to problem
(\ref{1.4})
\begin{equation}\label{5.1}
\aligned
J_t(u)=\frac12\int_{\Omega}
a(x)|\nabla u|^2-\int_{\Omega}a(x)k(x)|x-q|^{2\alpha}e^{-t\phi_1}e^u,
\,\quad\,\, u\in H_0^1(\Omega).
\endaligned
\end{equation}
For any points $\xi=(\xi_1,\ldots,\xi_m)\in\mathcal{O}_t(q)$, we
introduce the reduced energy
$\mathcal{F}_t:\mathcal{O}_t(q)\rightarrow\mathbb{R}$  by
\begin{equation}\label{5.2}
\aligned
\mathcal{F}_t(\xi)=J_t\big(u(\xi)
\big)=J_t\big(U(\xi)+\tilde{\phi}(\xi)
\big),
\endaligned
\end{equation}
where $U(\xi)$ is the approximation defined in (\ref{2.8})
and $\tilde{\phi}(\xi)(x)=\phi(\frac{\xi}{\varepsilon_0},\frac{x}{\varepsilon_0})$,
$x\in\Omega$, with $\phi=\phi_{\xi'}$  the unique solution
to problem (\ref{4.1}) given by Proposition 4.1. Define
\begin{equation*}\label{5.3}
\aligned
\mathcal{M}_m^t=\max\limits_{\xi=(\xi_1,\ldots,\xi_m)\in\mathcal{O}_t(q)}
\mathcal{F}_t(\xi).
\endaligned
\end{equation*}
Clearly,
for any $t$ large enough, the map
$\mathcal{F}_t(\xi)$ is differential in $\xi$ and hence the maximization problem
has a solution over $\mathcal{O}_t(q)$.

\vspace{1mm}
\vspace{1mm}
\vspace{1mm}
\vspace{1mm}

\noindent{\bf Proposition 5.1.}\,\,\,{\it
The maximization problem
\begin{eqnarray}\label{5.4}
\max\limits_{(\xi_1,\ldots,\xi_m)\in\mathcal{O}_t(q)}
\mathcal{F}_t(\xi_1,\ldots,\xi_m)
\end{eqnarray}
has a solution $\xi_t=(\xi_{1,t},\ldots,\xi_{m,t})\in\mathcal{O}_t^o(q)$, i.e., the interior of $\mathcal{O}_t(q)$.
}

\begin{proof}{\bf Step 1:}
With the choices  for  the parameters $\mu_0$ and $\mu_i$,
$i=1,\ldots,m$,
respectively given by {\upshape(\ref{2.12})} and {\upshape(\ref{2.13})},
let us claim that the following expansion   holds
\begin{eqnarray}\label{5.5}
J_t\big(U(\xi)\big)=
8\pi(1+\alpha)ta(q)
+8\pi t\sum_{i=1}^ma(\xi_i)
\phi_1(\xi_i)
+16\pi(2+\alpha)\sum_{i=1}^ma(\xi_i)\log|\xi_i-q|
+16\pi\sum_{i\neq j}^m
a(\xi_i)\log|\xi_i-\xi_j|+O\left(1\right)
\end{eqnarray}
uniformly for all points $\xi=(\xi_1,\ldots,\xi_m)\in\mathcal{O}_t(q)$ and for all $t$ large enough.
In fact,  observe  that by (\ref{2.8}) and (\ref{2.9}),
\begin{eqnarray}\label{5.6}
\frac12\int_{\Omega}a(x)|\nabla U|^2
=\frac12\int_{\Omega} \left[-\nabla\big(a(x)\nabla U_0\big)\right]U_0
+\sum_{j=1}^m\int_{\Omega} \left[-\nabla\big(a(x)\nabla U_0\big)\right]U_j
+\frac12\sum_{i,j=1}^m\int_{\Omega} \left[-\nabla\big(a(x)\nabla U_i\big)\right]U_j
&&\nonumber\\
=-\frac12\int_{\Omega}a(x)\big(u_0+H_0\big)\Delta u_0
-\sum_{j=1}^m\int_{\Omega}a(x)\big(u_j+H_j\big)\Delta u_0
-\frac12\sum_{i,j=1}^m\int_{\Omega}a(x)\big(u_j+H_j\big)\Delta u_i.
&&
\end{eqnarray}
Let us analyze the behavior of the first term.
By (\ref{2.4}), (\ref{2.5}) and  (\ref{2.10}) we get
$$
\aligned
-\int_{\Omega}a(x)\big(u_0+H_0\big)\Delta u_0
=\int_{\Omega}\frac{\,8\varepsilon_0^2\mu_0^2(1+\alpha)^2a(x)|x-q|^{2\alpha}\,}{\,(\varepsilon_{0}^2\mu_{0}^2+|x-q|^{2(1+\alpha)})^2\,}
\Big[\log
\frac{1}{(\varepsilon_{0}^2\mu_{0}^2+|x-q|^{2(1+\alpha)})^2}+
(1+\alpha)H(x,q)+O\left(\rho_0^\sigma v_0^\sigma\right)
\Big].
\endaligned
$$
Making the change of
variables $\rho_0v_0z=x-q$, we can derive that
$$
\aligned
-\int_{\Omega}a(x)\big(u_0+H_0\big)\Delta u_0
=\int_{\Omega_{\rho_0v_0}}\frac{8(1+\alpha)^2a(q+\rho_0v_0z)|z|^{2\alpha}}{(1+|z|^{2(1+\alpha)})^2}
\left[\log\frac{(\varepsilon_0\mu_0)^{-4}}{(1+|z|^{2(1+\alpha)})^2}
+(1+\alpha)H(q,q)+O\left(\rho_0^\sigma v_0^\sigma|z|^\sigma+\rho_0^\sigma v_0^\sigma\right)
\right],
\endaligned
$$
where $\Omega_{\rho_0v_0}=\frac{1}{\rho_0v_0}(\Omega-\{q\})$.  But
$$
\aligned
\int_{\Omega_{\rho_0v_0}}a(q+\rho_0v_0z)\frac{8(1+\alpha)^2|z|^{2\alpha}}{(1+|z|^{2(1+\alpha)})^2}
=8\pi(1+\alpha)a(q)+O\left(\rho_0v_0\right),
\endaligned
$$
and
$$
\aligned
\int_{\Omega_{\rho_0v_0}}a(q+\rho_0v_0z)\frac{8(1+\alpha)^2|z|^{2\alpha}}{(1+|z|^{2(1+\alpha)})^2}
\log\frac{1}{(1+|z|^{2(1+\alpha)})^2}=-16\pi(1+\alpha)a(q)+O(\rho_0v_0).
\endaligned
$$
Then
\begin{equation}\label{5.7}
\aligned
-\int_{\Omega}a(x)\big(u_0+H_0\big)\Delta u_0
=8\pi(1+\alpha)a(q)\big[(1+\alpha)H(q,q)-2-4\log(\varepsilon_0\mu_0)\big]
+O\left(\rho_0^\sigma v_0^\sigma\right).
\endaligned
\end{equation}
For the second term of (\ref{5.6}), by (\ref{2.2}), (\ref{2.4}), (\ref{2.5}), (\ref{2.11})
and  the change of
variables $\rho_0v_0z=x-q$ we get that for any $j=1,\ldots,m$,
\begin{eqnarray}\label{5.8}
-\int_{\Omega}a(x)\big(u_j+H_j\big)\Delta u_0
=\int_{\Omega}\frac{8\varepsilon_0^2\mu_0^2(1+\alpha)^2a(x)|x-q|^{2\alpha}}{\,(\varepsilon_{0}^2\mu_{0}^2+|x-q|^{2(1+\alpha)})^2\,}
\left[\log
\frac{1}{(\varepsilon_{j}^2\mu_{j}^2+|x-\xi_j|^2)^2}
+H(x,\xi_j)+O\left(\varepsilon_j^\sigma\mu_j^\sigma\right)
\right]dx
&&\nonumber\\
=\int_{\Omega_{\rho_0v_0}}\frac{8(1+\alpha)^2a(q+\rho_0v_0z)|z|^{2\alpha}}{(1+|z|^{2(1+\alpha)})^2}
\left[\log
\frac{1}{|q-\xi_j|^4}
+H(q,\xi_j)+O\left(\rho_0^\sigma v_0^\sigma|z|^\sigma\right)
\right.
\qquad\quad\,\,\,\,
&&\nonumber\\[1mm]
+O\left(\rho_0v_0t^\beta\big|z\big|\right)
+O\left(\varepsilon_j^\sigma\mu_j^\sigma\right)
\Big]dz
\qquad\qquad\qquad\qquad
\qquad\qquad\qquad\qquad
\qquad\qquad\quad\,\,\,\,
&&\nonumber\\[1mm]
=8\pi(1+\alpha)a(q)G(q,\xi_j)
+O\left(\rho_0^\sigma v_0^\sigma\right)
+O\left(\varepsilon_j^\sigma\mu_j^\sigma\right).
\qquad\qquad\qquad\qquad\qquad\quad\,\,\,\,
\qquad\qquad
&&
\end{eqnarray}
As for the last term of (\ref{5.6}), by (\ref{2.4}), (\ref{2.6}), (\ref{2.11})
and  the change of
variables $\varepsilon_i\mu_iz=x-\xi_i$ we observe that for any $i,\,j=1,\ldots,m$,
$$
\aligned
-\int_{\Omega}a(x)\big(u_j+H_j\big)\Delta u_i
=&\int_{\Omega}\frac{8\varepsilon_i^2\mu_i^2a(x)}{(\varepsilon_{i}^2\mu_{i}^2+|x-\xi_i|^2)^2}
\left[\log\frac{1}{(\varepsilon_{j}^2\mu_{j}^2+|x-\xi_j|^2)^2}+
H(x,\xi_j)+O(\varepsilon_j^\sigma\mu_j^\sigma)
\right]dx\\
=&\int_{\Omega_{\varepsilon_i\mu_i}}\frac{8a(\xi_i+\varepsilon_i\mu_iz)}{(1+|z|^2)^2}
\left[\log\frac{1}{(\varepsilon_{j}^2\mu_{j}^2+|\xi_i-\xi_j+\varepsilon_i\mu_iz|^2)^2}
+H(\xi_i,\xi_j)+
O\left(\varepsilon_i^\sigma\mu_i^\sigma|z|^\sigma\right)
+O\left(\varepsilon_j^\sigma\mu_j^\sigma\right)
\right]dz,
\endaligned
$$
where $\Omega_{\varepsilon_i\mu_i}=\frac{1}{\varepsilon_i\mu_i}(\Omega-\{\xi_i\})$.
Note that
$$
\aligned
\int_{\Omega_{\varepsilon_i\mu_i}}a(\xi_i+\varepsilon_i\mu_iz)\frac{8}{(1+|z|^2)^2}
=8\pi a(\xi_i)+O(\varepsilon_i\mu_i),
\endaligned
$$
and
$$
\aligned
\int_{\Omega_{\varepsilon_i\mu_i}}a(\xi_i+\varepsilon_i\mu_iz)\frac{8}{(1+|z|^2)^2}
\log\frac{1}{(1+|z|^2)^2}=-16\pi a(\xi_i)+O(\varepsilon_i\mu_i).
\endaligned
$$
Then for all $i,\,j=1,\ldots,m$,
\begin{equation}\label{5.9}
\aligned
-\int_{\Omega}a(x)\big(u_j+H_j\big)\Delta u_i
=\left\{\aligned
&8\pi a(\xi_i) \big[H(\xi_i,\xi_i)-2-4\log(\varepsilon_i\mu_i)\big]
+O\left(\varepsilon_i^\sigma\mu_i^\sigma\right)
\quad\,\,\forall\,\,\,i=j,\\[2.5mm]
&8\pi
a(\xi_i) G(\xi_i,\xi_j)+O\left(\varepsilon_i^\sigma\mu_i^\sigma+\varepsilon_j^\sigma\mu_j^\sigma\right)
\qquad\qquad\qquad\,
\forall\,\,\,i\neq j.
\endaligned\right.
\endaligned
\end{equation}
On the other hand, by (\ref{2.18}), (\ref{2.19}),  (\ref{2.20}) and  the change of variables
$x=\varepsilon_0y=e^{-\frac12t}y$, we obtain
$$
\aligned
\int_{\Omega}a(x)k(x)|x-q|^{2\alpha}e^{-t\phi_1}e^U=&
\int_{\Omega_t}a(\varepsilon_0 y)k(\varepsilon_0 y)|\varepsilon_0 y-q|^{2\alpha}e^{-t\big[\phi_1(\varepsilon_0 y)-1\big]}e^{U(\varepsilon_0 y)-2t}dy
\\[1mm]
=&\left(\int_{\Omega_t\setminus\big[\bigcup_{i=1}^mB_{\frac1{\varepsilon_0 t^{2\beta}}}(\xi'_i)\cup
B_{\frac1{\varepsilon_0 t^{2\beta}}}(q')\big]}
+\int_{B_{\frac1{\varepsilon_0 t^{2\beta}}}(q')}
+\sum\limits_{i=1}^m\int_{B_{\frac1{\varepsilon_0 t^{2\beta}}}(\xi'_i)}\right)a(\varepsilon_0 y)W dy.
\endaligned
$$
By (\ref{2.30}),  (\ref{2.28})  and (\ref{2.32}) we obtain
$$
\aligned
\int_{\Omega_t\setminus\big[\bigcup_{i=1}^mB_{\frac1{\varepsilon_0 t^{2\beta}}}(\xi'_i)\cup
B_{\frac1{\varepsilon_0 t^{2\beta}}}(q')\big]}
a(\varepsilon_0 y)W dy=&
\int_{\Omega_t\setminus\big[\bigcup_{i=1}^mB_{\frac1{\varepsilon_0 t^{2\beta}}}(\xi'_i)\cup
B_{\frac1{\varepsilon_0 t^{2\beta}}}(q')\big]}
O\left(\frac{\varepsilon_0^2e^{-t\phi_1(\varepsilon_0 y)}}
{\,|\varepsilon_0 y-q|^{4+2\alpha}\,}\prod_{i=1}^m\frac{1}{\,|\varepsilon_0 y-\xi_i|^4\,}
\right)dy\\[2mm]
=&\,O\left(1\right),
\endaligned
$$
and
$$
\aligned
\int_{B_{\frac1{\varepsilon_0 t^{2\beta}}}(q')}a(\varepsilon_0 y)Wdy=&
\int_{B_{\frac1{\varepsilon_0 t^{2\beta}}}(q')}\left(\frac{\varepsilon_0}{\rho_0v_0}\right)^2
\frac{8(1+\alpha)^2a(\varepsilon_0 y)\big|\frac{\varepsilon_0 y-q}{\rho_0v_0}\big|^{2\alpha}}{\,\big(1+\big|\frac{\varepsilon_0 y-q}{\rho_0v_0}\big|^{2(1+\alpha)}\big)^2\,}
\Big[1+O\left(\varepsilon_0^\sigma|y-q'|^\sigma\right)
+
o\left(1\right)
\Big]dy\\[2mm]
=&\,8\pi(1+\alpha)a(q)+o(1),
\endaligned
$$
and for any $i=1,\ldots,m$,
$$
\aligned
\int_{B_{\frac1{\varepsilon_0 t^{2\beta}}}(\xi'_i)}a(\varepsilon_0 y)Wdy=
\int_{B_{\frac1{\varepsilon_0 t^{2\beta}}}(\xi'_i)}\frac{1}{\gamma_{i}^2}
\frac{8a(\varepsilon_0 y)}{\big(1+\big|\frac{y-\xi'_i}{\gamma_{i}}\big|^2\big)^2}
\Big[1+O\left(\varepsilon_0^\sigma|y-\xi'_i|^\sigma\right)
+
o\left(1\right)
\Big]dy
=8\pi a(\xi_i)+o(1).
\endaligned
$$
Then
\begin{equation}\label{5.10}
\aligned
\int_{\Omega}a(x)k(x)|x-q|^{2\alpha}e^{-t\phi_1}e^U=O\left(1\right).
\endaligned
\end{equation}
Hence  by (\ref{5.1}), (\ref{5.6})-(\ref{5.10}) we conclude that
$$
\aligned
J_t\left(U(\xi)\right)=-16\pi(1+\alpha)a(q)\log(\varepsilon_0\mu_0)-16\pi
\sum_{i=1}^ma(\xi_i)\log(\varepsilon_i\mu_i)
+8\pi(1+\alpha)\sum\limits_{i=1}^{m}a(q)G(q,\xi_i)
+4\pi\sum_{i\neq j}^ma(\xi_i)G(\xi_i,\xi_j)
+O(1),
\endaligned
$$
which, together with the definitions  of $\varepsilon_0$, $\varepsilon_i$ in (\ref{2.7})
and the choices of $\mu_0$, $\mu_i$ in  (\ref{2.12})-(\ref{2.13}),
implies that expansion (\ref{5.5}) holds.

\vspace{1mm}

{\bf Step 2:}
For any  $t$ large enough,  we claim that
the
following expansion holds
\begin{equation}\label{5.11}
\aligned
\mathcal{F}_t(\xi)=J_t\big(U(\xi)\big)+o(1)
\endaligned
\end{equation}
uniformly on points $\xi=(\xi_1,\ldots,\xi_m)\in\mathcal{O}_t(q)$.
Indeed, if we define
\begin{equation}\label{5.12}
\aligned
I_t(\omega)=\frac12\int_{\Omega_t}
a(\varepsilon_0 y)|\nabla \omega|^2-\int_{\Omega_t}|\varepsilon_0y-q|^{2\alpha}a(\varepsilon_0 y)\kappa(y,t)e^\omega,
\,\quad\,\,\omega\in H_0^1(\Omega_t),
\endaligned
\end{equation}
then by (\ref{2.7}) and (\ref{2.19}),
$$
\aligned
\mathcal{F}_t(\xi)-
J_t\big(U(\xi)\big)=I_t\big(V(\xi')+\phi_{\xi'}
\big)-I_t\big(V(\xi')
\big).
\endaligned
$$
Using  $DI_t(V+\phi_{\xi'})[\phi_{\xi'}]=0$, a Taylor expansion and an integration by parts  give
$$
\aligned
\mathcal{F}_t(\xi)-
J_t\big(U(\xi)\big)=&\int_0^1D^2I_t(V+\tau\phi_{\xi'})[\phi_{\xi'}]^2(1-\tau)d\tau
\\
=&\int_0^1\left\{
\int_{\Omega_t}a(\varepsilon_0 y)
\Big[\big(N(\phi_{\xi'})+E\big)\phi_{\xi'}+
W\big(1-e^{\tau\phi_{\xi'}}\big)\phi_{\xi'}^2\Big]
\right\}(1-\tau)d\tau.
\endaligned
$$
From the estimates in Lemma $3.4$ and Proposition $4.1$, we find
$$
\aligned
\mathcal{F}_t(\xi)-
J_t\left(U(\xi)\right)=O\left(t
\max\left\{\big(\rho_0v_0\big)^{\min\{2\sigma,4(\alpha-\hat{\alpha})\}},\,(\varepsilon_1\mu_1)^{2\sigma},\,\ldots,\,
(\varepsilon_m\mu_m)^{2\sigma},\,
\|e^{-t\phi_1}\|_{L^{\infty}(\Omega_t)}\right\}\right)
=o\left(1\right).
\endaligned
$$
The continuity in $\xi$ of this  expression is inherited from that of $\phi_{\xi'}$  in the $L^\infty$ norm.

\vspace{1mm}

{\bf Step 3:} Proof of Proposition 5.1.
Let $\xi_t=(\xi_{1,t},\ldots,\xi_{m,t})$
be the maximizer of $\mathcal{F}_t$ over $\mathcal{O}_t(q)$. We need to prove that $\xi_t$
belongs to the interior of $\mathcal{O}_t(q)$.
First, we obtain a lower bound for $\mathcal{F}_t$ over $\mathcal{O}_t(q)$.
Let us fix the point
$q$ as an isolated local maximum point of $a(x)\phi_1$ in $\Omega$ and set
$$
\aligned
\xi^0_i=q+
\frac{1}{\sqrt{t}}\widehat{\xi}_i,
\endaligned
$$
where $\widehat{\xi}=(\widehat{\xi}_1,\ldots,\widehat{\xi}_m)$ forms
an $m$-regular polygon in $\mathbb{R}^2$.
Obviously,
$\xi^0=(\xi^0_1,\ldots,\xi^0_m)\in\mathcal{O}_t(q)$
because  $\beta>1$ and $a(\xi^0_i)\phi_1(\xi^0_i)=a(q)\phi_1(q)+O(t^{-1})$.
Then
\begin{eqnarray}\label{5.13}
\max\limits_{\xi\in\mathcal{O}_t(q)}
\mathcal{F}_t(\xi)\geq
8\pi(1+\alpha)ta(q)
+8\pi t\sum_{i=1}^ma(\xi_i^0)
\phi_1(\xi_i^0)
+16\pi(2+\alpha)\sum_{i=1}^ma(\xi_i^0)\log|\xi_i^0-q|
+16\pi\sum_{i\neq j}^m
a(\xi_i^0)\log|\xi_i^0-\xi_j^0|+O\left(1\right)
&&\nonumber\\
\geq
8\pi(m+1+\alpha)a(q)t
-8\pi (m+1+\alpha))\big[a(\xi_1^0)+\cdots+a(\xi_m^0)\big]
\log t+O(1).
\qquad\qquad\qquad\qquad
\qquad\qquad\qquad\quad
&&
\end{eqnarray}
Next, we suppose $\xi_t=(\xi_{1,t},\ldots,\xi_{m,t})\in\partial\mathcal{O}_t(q)$.
Then there  exist four possibilities:\\
C1. \,There exists an $i_0$ such that
$\xi_{i_0,t}\in\partial B_d(q)$, in which case,
$a(\xi_{i_0,t})\phi_1(\xi_{i_0,t})\leq
a(q)\phi_1(q)-d_0$ for some
 $d_0>0$
 \indent\indent
 independent of $t$;\\
C2. \,There exists an $i_0$ such that
$a(\xi_{i_0,t})\phi_1(\xi_{i_0,t})=a(q)\phi_1(q)-\frac{1}{\sqrt{t}}$;\\
C3. \,There exist indices $i_0$, $j_0$,
$i_0\neq j_0$ such that
$|\xi_{i_0,t}-\xi_{j_0,t}|=t^{-\beta}$;\\
C4. \,There exists an  $i_0$
such that
$|\xi_{i_0,t}-q|=t^{-\beta}$.\\
For the first case,   we have
\begin{eqnarray}\label{5.14}
\max\limits_{\xi\in\mathcal{O}_t(q)}
\mathcal{F}_t(\xi)=\mathcal{F}_t(\xi_t)
\leq
8\pi(1+\alpha)ta(q)
+8\pi t\left[(m-1)a(q)
\phi_1(q)+a(q)
\phi_1(q)-d_0\right]+O\big(\log t\big)&&\nonumber\\[1mm]
=8\pi(m+1+\alpha)a(q)t-8\pi  d_0 t+O\big(\log t\big),
\qquad\qquad\qquad\qquad\qquad\quad\quad\,\,
&&
\end{eqnarray}
which contradicts to (\ref{5.13}).
This shows that $a(\xi_{i,t})\phi_1(\xi_{i,t})\rightarrow a(q)\phi_1(q)$. Using the assumption of $a(x)\phi_1$ over $\oo$,
we deduce $\xi_{i,t}\rightarrow q$ for all $i=1,\ldots,m$.

For the second case,   we have
\begin{eqnarray}\label{5.14}
\max\limits_{\xi\in\mathcal{O}_t(q)}
\mathcal{F}_t(\xi)=\mathcal{F}_t(\xi_t)
\leq
8\pi(1+\alpha)ta(q)
+8\pi t\left[(m-1)a(q)
\phi_1(q)+a(q)
\phi_1(q)-\frac1{\sqrt{t}}\right]+O\big(\log t\big)&&\nonumber\\[1mm]
=8\pi(m+1+\alpha)a(q)t-8\pi\sqrt{t}+O\big(\log t\big),
\qquad\qquad\qquad\qquad\qquad\qquad\quad\,\,\,
&&
\end{eqnarray}
which contradicts to (\ref{5.13}).
For the third case, we have
\begin{equation}\label{5.15}
\aligned
\max\limits_{\xi\in\mathcal{O}_t(q)}
\mathcal{F}_t(\xi)=\mathcal{F}_t(\xi_t)
\leq
8\pi(m+1+\alpha) ta(q) -16\pi\beta
a(\xi_{i_0,t})\log t
+O\big(1\big).
\endaligned
\end{equation}
For the last case, we have
\begin{equation}\label{5.16}
\aligned
\max\limits_{\xi\in\mathcal{O}_t(q)}
\mathcal{F}_t(\xi)=\mathcal{F}_t(\xi_t)
\leq
8\pi(m+1+\alpha) t a(q) -16\pi(2+\alpha)\beta
a(\xi_{i_0,t})\log t
+O\big(1\big).
\endaligned
\end{equation}
Combining (\ref{5.15})-(\ref{5.16}) with (\ref{5.13}), we give
\begin{equation}\label{5.17}
\aligned
16\pi(2+\alpha)\beta a(\xi_{i_0,t})\log t
+O\big(1\big)\leq
8\pi (m+1+\alpha))\big[a(\xi_1^0)+\cdots+a(\xi_m^0)\big]
\log t+O(1),
\endaligned
\end{equation}
which  is impossible by the choice of $\beta$ in (\ref{2.3}).
\end{proof}

\section{Proof of Theorem 1.1}
\noindent {\bf Proof of Theorem 1.1.}
From Proposition 4.1 it follows that for
any points
$\xi=(\xi_1,\ldots,\xi_m)\in\mathcal{O}_t(q)$ and any $t$ large enough,
there is a function
$\phi_{\xi'}$, which solves
$$
\aligned
-\Delta_{a(\varepsilon_0 y)}\big(V(\xi')+\phi_{\xi'}\big)-|\varepsilon_0y-q|^{2\alpha}\kappa(y,t)e^{V(\xi')+\phi_{\xi'}}
=\frac{1}{a(\varepsilon_0 y)}\sum\limits_{i=1}^m\sum\limits_{j=1}^2c_{ij}(\xi')\chi_iZ_{ij},
\qquad
\int_{\Omega_t}\chi_iZ_{ij}\phi_{\xi'}=0
\,\,\,\,\,
\endaligned
$$
for some coefficients $c_{ij}(\xi')$,
$i=1,\ldots,m$, $j=1,2$.
So, in order to
find a solution to  equation
(\ref{2.17}) and hence to the original problem  (\ref{1.4}),
we need to match  $\xi'$ with
the coefficients $c_{ij}(\xi')$  such that
\begin{equation}\label{6.1}
\aligned
c_{ij}(\xi')=0\,\,\quad\,\,\textrm{for all}\,\,\,i=1,\ldots,m,\,\,j=1,2.
\endaligned
\end{equation}
According to Proposition 5.1, there exits
a $\xi_t=(\xi_{1,t},\ldots,\xi_{m,t})\in\mathcal{O}_t^o(q)$
that achieves the maximum for the maximization problem (\ref{5.4}).
Let
$\omega_t=V(\xi'_t)+\phi_{\xi'_t}$.
Then we have
\begin{equation}\label{6.2}
\aligned
\partial_{\xi_{kl}}\mathcal{F}_t(\xi_t)=0\,\,\quad\,\,\textrm{for all}\,\,\,k=1,\ldots,m,\,\,l=1, 2.
\endaligned
\end{equation}
From (\ref{5.1}), (\ref{5.2}) and (\ref{5.12}) we get
$$
\aligned
\partial_{\xi_{kl}}\mathcal{F}_t(\xi_t)=&\,\partial_{\xi_{kl}}J_t\big(U(\xi_t)+\tilde{\phi}(\xi_t)
\big)
=\,\frac1{\varepsilon_0}
\partial_{\xi'_{kl}}I_t\big(V(\xi_t')+\phi_{\xi_t'}
\big)\\
=&\,\frac1{\varepsilon_0}\left\{
\int_{\Omega_t}a(\varepsilon_0 y)\nabla\omega_t\nabla\left[\partial_{\xi'_{kl}}V(\xi'_t)+\partial_{\xi'_{kl}}\phi_{\xi'_t}
\right]-\int_{\Omega_t}a(\varepsilon_0 y)|\varepsilon_0y-q|^{2\alpha}\kappa(y,t)e^{\omega_t}\left[\partial_{\xi'_{kl}}V(\xi'_t)+\partial_{\xi'_{kl}}\phi_{\xi'_t}
\right]
\right\}.
\endaligned
$$
Then for all $k=1,\ldots,m$
and $l=1, 2$,
$$
\aligned
\sum\limits_{i=1}^m\sum\limits_{j=1}^2c_{ij}(\xi'_t)\int_{\Omega_t}\chi_iZ_{ij}
\left[\partial_{\xi'_{kl}}V(\xi'_t)+\partial_{\xi'_{kl}}\phi_{\xi'_t}
\right]=0.
\endaligned
$$
Notice  that
$\partial_{\xi'_{kl}}V(\xi'_t)=Z_{kl}+o\big(1\big)$
and $\partial_{\xi'_{kl}}\phi_{\xi'_t}
=o(1)$ with $o(1)$ sufficiently small in the sense of the $L^{\infty}$ norm
as $t\rightarrow+\infty$. Therefore, we find  that $D_{\xi}\mathcal{F}_t(\xi)=0$
implies  the validity of a system of equations of the form
$$
\aligned
\sum\limits_{i=1}^m\sum\limits_{j=1}^2c_{ij}(\xi_t')\int_{\Omega_t}\chi_iZ_{ij}\big[
Z_{kl}(y)+o\left(1\right)
\big]=0,\,\ \ \ \ \,\,\,\,\,k=1,\ldots,m,\,\,l=1,2.
\endaligned
$$
Clearly, the  matrix of this  system is  diagonally dominant and then
$c_{ij}(\xi_t')=0$ for all $i=1,\ldots,m$, $j=1,2$.
As a result, we find a solution $u_t$ of problem (\ref{1.4}) in the form
$U(\xi_t)+\tilde{\phi}(\xi_t)$ with
the qualitative properties  as predicted in
Theorem 1.1.  \,\,\,\,$\square$

\end{document}